\documentclass[11pt,a4paper]{article}
\usepackage{geometry}
\usepackage{amsmath,amssymb,amsfonts,amsthm}
\usepackage{mathtools}
\usepackage[colorlinks=true,linkcolor=magenta,citecolor=teal]{hyperref}
\usepackage{cleveref}
\crefname{appsec}{Appendix}{Appendices}
\usepackage{multirow}
\usepackage{algorithm}
\usepackage{algpseudocode}
\usepackage[table]{xcolor}
\usepackage{booktabs}
\usepackage{rotating, graphicx}
\usepackage{lscape}
\usepackage{longtable}
\usepackage{natbib}
\usepackage{appendix}
\usepackage{bm}
\usepackage{bbold}
\usepackage{algorithmicx}
\usepackage{tikz}
\usepackage{caption,subcaption}

\allowdisplaybreaks
\definecolor{LG}{gray}{0.9}
\definecolor{Gray}{gray}{0.85}

\newcommand{\mtc}{\mathcal}
\newcommand{\pr}{\mathsf}
\newcommand{\rv}{\boldsymbol}

\newtheorem{proposition}{Proposition}
\newtheorem{example}{Example}

\title{Solution of Stochastic Facility Location Problems with Combinatorially many Decision-Dependent Distributions}
\author{Giovanni Pantuso \\
Department of Mathematical Sciences \\
University of Copenhagen \\
Universitetsparken 5, 2100, Copenhagen, Denmark\\
\href{mailto:gp@math.ku.dk}{gp@math.ku.dk}}
\date{}

\begin{document}

\maketitle

\begin{abstract}
  This article describes a model and an exact solution method for facility location problems with decision-dependent uncertainties.
  The model allows characterizing the probability distribution of the random elements as a function of the choice of open facilities.
  This, in turn, generates a combinatorial number of potential distributions of the random elements.
  Though general in the relationship between location decisions and distributions, the proposed model is, however, exponential in size.
  We show that the problem can be solved efficiently by a recent finitely convergent method for stochastic programs with decision-dependent uncertainty,
  for which we tight prove cutting planes and efficient valid inequalities. 
  Extensive tests show that facility location problems with up to $2^{17}$ potential distributions and hundreds of thousand scenarios are solved within minutes.
  These results indentify a promising solution strategy for other combinatorial optimization problems characterized by decision-dependent uncertanty. 
  {\bf Keywords:} Stochastic Programming, Endogenous Uncertainty, Facility Location.
\end{abstract}

\section{Introduction}
Discrete facility location problems are extensively explored in optimization studies.
The research literature is rich of key methodological developments, see e.g., \cite{KueH63,Man64,GeoM78,Van86,HolRY99,SahS07,FarSA10,BasSG15,LiuSWM22} and applications, see e.g., \
\cite{OweD98,MelNS09,LiZZW11,AhmSS17,TerAF12,MarS94}.
The books \citet{Das13} and \cite{LapNS19} provide a comprehensive overview of the subject.

Location decisions are often made in the presence of uncertainties.
The information guiding these decisions is seldom available when decisions are made.
This information includes factors such as demand levels, transportation costs, customer presence and location, and commodity prices.
The literature widely recognizes that uncertainty may affect this information in various ways.
The vast body of literature on facility location problems under uncertainty is thoroughly summarized in \cite{Lou86}, \cite{Sny06}, \cite{CorS19} and \cite{SalW24}. 
In the existing studies the uncertainty is almost exclusively exogenous, that is, independent of location decisions.
This setup represents many, but perhaps not all, facility location problems.
In this paper we address situations where the uncertainty is decision-dependent, i.e., endogenously determined by location decisions.

Location decisions may subtly influence the uncertainties.
For instance, the effect of proximity on demand has been widely reconized, starting from the classical Hotelling's model \cite{Hot29}, see e.g., \cite{WanG17,ShrB22,JekBE01,PauMJWTPG10,HoP95}.
As an example, \cite{ShrB22} estimate that a $10$\% reduction in retail store distance leads to $1.9$\% higher expenditure.
These evidences are not limited to the retail domain, see e.g., \cite{JekBE01,EvaH05,HarMEMPS19,LevLSP11,Pra09,AlwD97}.
In a similar way, the distribution of transportation costs may be influenced by the distance to customers or suppliers,
the transport infrastrucuture and economies of scale. Geopolitical factors may have an impact on volatity. 
An exhaustive investigation of the intricate relation between location decisions and various economic inputs is beyond the scope of this paper.
In what follows, we focus on optimization techniques that allow solving the resulting problems.

The literature on facility location problems with decision-dependent uncertainty counts, to our knowledge, only the recent article by \cite{BasAS21}.
The authors study a facility location problem assuming that demand is endogenously influenced by location decisions.
Particularly, the mean and variance of the demand distributions are functions of facility location decisions.
They propose a distributionally-robust optimization model to address mispecifications in the demand distribution.
The objective is to minimize location cost and the worst-case expected cost of transportation and unsatisfied demand.
Location decisions determine an ambiguity set which contains the potential demand distributions that match the first two moments.
The resulting model optimizes against the worst-case cost among all the possible distributions in the resulting ambiguity set.
The authors present an exact MILP reformulation for a special case where demand mean and variance are piecewise linear functions of location decisions.
In an extensive computational study the authors demonstrate the superiority of the proposed model to traditional robust optimization and stochastic programs in terms of increase of profit and demand satisfied.

Stochastic programs with decision-dependent uncertainty are notoriously challenging.
The literature offers a number of practical applications, see e.g., \cite{Ahm00,VisSF04,Fla10,PeeSGV10,LauPK14,TonFR12,EscGMU18}. However, general-purpose methodology is sparse. 
\cite{Pfl90} provide a convergent algorithm for a Markovian random process whose transition probabilities depend on the decision variables of the optimization problem.
\cite{JonWW98} consider the case where decisions influence both the probability measure and the timing of the observation. 
The authors assume that the set of probability measures which can be enforced by decisions is finite and countable.
This is also the case in our paper.
The authors show that the problem can be solved by choosing the best among the stochastic programs determined by a choice of a distribution and propose an implicit enumeration algorithm.
\cite{GoeG06} consider multi-stage stochastic programs where decisions have an impact on the resolution time of the underlying stochastic process.
The authors provide disjunctive programming formulations and a provide a Lagrangian relaxation-based branch \& bound framework for its solution. 
\cite{Pan21} also considers multi-stage problems and extends the model of \cite{JonWW98} by allowing decisions at all stages to determine the probability measure for the later stages. The author proposes a formulation that avoids non-anticipativity constraints, see e.g., \cite{ApaG16,HooM16,MooM18} and an algorithm for a special case.
\cite{HelBT18} discuss several ways of modeling the interplay between the decision variables and the parameters of the underlying probability distributions.
Examples are cases where prior probabilities are distorted through affine transformations, or combined using convex combinations of several probability distributions.
Finally, \citet{PanH25} propose an extension of the L-Shaped method for two-stage stochastic programs with decision-dependent uncertainty.
The method is based on identifying the enforced distribution for each candidate first-stage solution prior to separate decision-dependent cuts.
It relies on the assumption that there exists finitely many probability distributions that may be enforced.
This assumption naturally holds for many combinatorial optimization problems with decision-dependent uncertainty.
In fact, the set of decisions is often countable and finite, though large. Hence, each (subset of) decision(s) may enforce a different distribution of the data.
In particular, as we will show, this assumption holds in the context of location decisions.
The consequence is, however, that the number of potentially distributions is combinatorially large.

In this article we propose a general model for two-stage facility location problems with decision-dependent uncertainty.
The model allows specifying probability distributions of the second-stage parameters that are fully characterized by location decisions.
We show that this problem gives rise to countably (though exponentially) many potential distributions.
To solve the problem we propose an exact decomposition algorithm based on the recent method described in \citet{PanH25}.
In particular, we derive distribution-specific cuts for which we show validity.
We also show that the proposed cuts ensure finite convergence of the method.
We report on extensive numerical tests on a two-stage facility location problem where demand is stochastic and depends on the choice of open facilities.
The framework is however general and can immediately be applied to facility location problems characterized by any type of relationship between decisions and uncertainties.
To the best of our knowledge the proposed method is the first general-purpose method for two-stage facility location problems with decision-dependent uncertainty.

This article extends the available literature in a number of ways.
First, we extend the available literature with a study on facility location decisions with decision-dependent uncertainties.
This adds to the pioneer work of \cite{BasAS21} where location decisions influence the distribution of demand.
Second, we allow the distribution of any (all) of random parameter(s) to be determined by location decisions.
Third, we allow the probability distributions to be fully characterized by location decision. Hence, we do not restrict ourselves to a particular type of dependence.
Fourth, in addition to a general monolithic formulation, we provide an exact, general-purpose, decomposition algorithm. 
Fifth, we report on extensive numerical experiments. These show that the monolithic formulation can handle only relatively small instances.
On the other hand the proposed method is able to solve to provable optimality significantly larger instances. Hence, we significally extend the range of practically tractable problems.
Finally, we illustrate a promising solution strategy for combinatorial optimization problems characterized by decision-dependent uncertainty.

The remainder of this article is organized as follows.
In \Cref{sec:model} we introduce a rather abstract and general mathematical model for facility location decisions under endogenous uncertainty.
In \Cref{sec:method} we describe a general-purpose exact solution method.
In \Cref{sec:exp} we report on extensive numerical tests on a facility location problem with endogenous stochastic demand.
In \Cref{sec:ex} we discuss the potential impact of accounting for decision-dependent uncertainty.
Finally, we draw conclusions in \Cref{sec:conclusions}.

\section{General Model}\label{sec:model}

In this section we introduce the problem in an abstract general way.
 We assume a two-stage decision process. Facilities are open in the first decision stage.
At the time location decisions are made, several parameters of the problem (e.g., costs and demands) are uncertain.
Facilities are operated in the second decision stage.
At that time, facility locations decisions have been made and implemented and uncertain problem parameters have materialized.

We will allow the distribution of the uncertainty parameters to depend on which subsets of facilities have been opened.
This choice is without loss of generality. Subsets may, in fact, contain only one facility.
Hence, this framework allows modeling problems where each individual selection of open facilities determines a different distribution.

Let $\mtc{I}$ be the set of potential facility locations. Opening a facility at location $i$ determines a fixed cost $F_i$.
The effect of location decisions depends on a number of parameters which will fully materialize only after location decisions have been made, i.e., in the second decision stage.  Examples of these are transportation costs between open facilities and customers, customer demands and effective capacities at facilities or links.
We model such parameters as a (multidimensional) random variable $\rv{\xi}:(\Omega,\mtc{F})\to(\mathbb{R}^{N},\mtc{B}^N)$ from a measurable space $(\Omega,\mtc{F})$ to $\mathbb{R}^N$ equipped with the Borel $\sigma$-algebra $\mtc{B}^N$. 
The probability distribution of $\rv{\xi}$ on $\mathbb{R}^N$ is, however, fully specified by the location decisions made, as we explain next.

Let $x:=(x_i)_{i\in\mtc{I}}\in\{0,1\}^{|\mtc{I}|}$ represent facility location decisions.
To allow for a rather general setup -- on whose practical relevance we comment later --
we assume that there exists a partition $(\mtc{I}_z)_{z\in\mtc{Z}}\subset 2^{\mtc{I}}$ of $\mtc{I}$,
with $\mtc{Z}$ being a finite set of indices of subsets, $|\mtc{Z}|\leq|\mtc{I}|$, and $2^{\mtc{I}}$ being the power set of $\mtc{I}$.
Note that $|\mtc{Z}|\leq|\mtc{I}|$ entails that the number of subsets in the partition may in principle be equal to the number of facilities.
For each $z\in\mtc{Z}$ we define an indicator function:
$$\mathbb{1}_{\mtc{I}_z}(x):=\begin{cases}1, \text{ iff } \sum_{i\in\mtc{I}_z}x_i>0\\
  0, \text{otherwise}
\end{cases}$$
Thus, $\mathbb{1}_{\mtc{I}_z}(x)$ indicates whether decision $x$ opens any facility in subset $\mtc{I}_z$ or not.
Observe, that $\mtc{Z}$ can represent a variety of configurations. These range from $|\mtc{Z}|=1$, where $\mtc{I}$ is essentially not partitioned, to $|\mtc{Z}|=|\mtc{I}|$ where each location defines a subset.
We anticipate that this, in turn, allows us to make the distribution of $\rv{\xi}$ dependent on the sets with open facilities (to already gain some practical intuition one may think of $\mtc{Z}$ as a set of geographical zones or facility types).
However, by acting on $\mtc{Z}$ one can model this dependence at any level of granularity.
This includes the special (and most general) case $|\mtc{Z}|=|\mtc{I}|$ where each facility is unrelated (say, e.g., geographically or otherwise) to the others and
any choice of open facilities determines a different distribution. In any intermediate case we have instead that any choice of ``open subsets'' determines a distribution.
 
Particularly, we let the distribution of $\rv{\xi}$ be determined by the $x$ throrugh $(\mathbb{1}_{\mtc{I}_z}(x))_{z\in\mtc{Z}}\in\{0,1\}^{\vert\mtc{Z}\vert}$, that is, by the open subsets.
Since the number of possible occurrences of $(\mathbb{1}_{\mtc{I}_z}(x))_{z\in\mtc{Z}}$ is finite, precisely $2^{|\mtc{Z}|}$, there can only be as many distributions of $\rv{\xi}$.
Let $\mtc{D}$ be the set of possible distributions.
More formally, one may say that there exists a set $\mtc{D}$ of probability measures $\mu_d$ on $(\Omega,\mtc{F})$ which induce distributions $\pr{P}_d$ on $\mathbb{R}^N$ for $\rv{\xi}$, hence making subsets of the possible realizations of $\rv{\xi}$ more or less likely, and potentially assigning probability zero.
Thus, some realizations of $\rv{\xi}$ may occur in some distributions but have probability zero in other.
This level of detail will however not be necessary.
Each distribution $d\in\mtc{D}$ is enforced by the activation of a unique collection of subsets $\mtc{I}_z$ for a distribution-specific collection $\mtc{Z}_d$ of indices.
That is, distribution $\mathsf{P}_d$ for $d\in\mtc{D}$ is enforced by $x$ if and only if $\sum_{z\in\mtc{Z}_d}\mathbb{1}_{\mtc{I}_z}(x)-\sum_{z\in\mtc{Z}\setminus\mtc{Z}_d}\mathbb{1}_{\mtc{I}_z}(x)=|\mtc{Z}_d|$, thus iff only facilities in the subsets in $\mtc{Z}_d$ are open and the remaining facilities are closed.
Hence, we define 
$$\mathbb{1}_{d}(x):=\begin{cases}1, \text{ iff } \sum_{z\in\mtc{Z}_d}\mathbb{1}_{\mtc{I}_z}(x)-\sum_{z\in\mtc{Z}\setminus\mtc{Z}_d}\mathbb{1}_{\mtc{I}_z}(x)=|\mtc{Z}_d|\\
  0, \text{otherwise}
\end{cases}$$

The specific ways in which active subsets in $\mtc{Z}$ determine conditional probability distributions $\mathsf{P}_d$, $d\in\mtc{D}$, for $\rv{\xi}$ is clearly application-specific. 
In a rather general setting one may consider for example a probability distribution with parameters in $\mathbb{R}^p$ determined by the specific activation of subsets.
Then, we can let $\theta:\{0,1\}^{\vert \mtc{Z}\vert}\to\mathbb{R}^p$ define the vector of parameters of a distribution $\mathsf{P}_\theta$ as a function of the selection of subsets.
Hence, we obtain distributions $\mathsf{P}_{\theta\left(\mathbb{1}_{\mtc{I}_{z_1}}(x),\ldots,\mathbb{1}_{\mtc{I}_{\vert\mtc{Z}\vert}}(x)\right)}$ for every choice of $x$, where the parameters are specifically determined by the subsets of $\mtc{I}$ activated by $x$.
Since the selection of active subsets is finite, we may enumerate distributions $\mathsf{P}_d$ as follows.
Define, for all $d\in\mtc{D}$, $\mathbb{1}_{\mtc{Z}_d}:\mtc{Z}\to \{0,1\}$ as 
$$\mathbb{1}_{\mtc{Z}_d}(z):=
\begin{cases}
  1, \text{ iff } z\in\mtc{Z}_d\\
  0, \text{ otherwise }
\end{cases}$$
Then we obtain, for all $d\in\mtc{D}$
$$\mathsf{P}_d=\mathsf{P}_{\theta\left(\mathbb{1}_{\mtc{Z}_d(z_1)},\ldots,\mathbb{1}_{\mtc{Z}_d(z_{\vert\mtc{Z}\vert})}\right)}$$
A concrete example is provided in \Cref{sec:exp}. In an even more general framework, one may assume that $\mathsf{P}_d$ takes a different functional form for each $d\in\mtc{D}$.

Finally, let $Q(x,\xi)$ represent the second-stage (e.g., operating) costs of a choice of open facilities $x$, given a realization $\xi$ of the random variable $\rv{\xi}$.
We leave the specification of $Q(x,\xi)$ open. In general, we assume $Q(x,\xi)$ is the result of an optimization problem. Therefore, it is not necessarily available in closed form. For example, $Q(x,\xi)$ determines the optimal allocation of customers to the open facilities decided by $x$, given a realization $\xi$ of the random parameters. A possible specification is provided in \Cref{sec:exp}. 
Given this notation we can now formulate a general \textit{facility location problem under decision-dependent uncertainty} as follows:

  \begin{align}\label{eq:flp:compact}
    \min_{x\in\{0,1\}^{|\mtc{I}|}}  &F^\top x+\sum_{d\in\mtc{D}}\mathbb{1}_d(x)\mathbb{E}_{\mathsf{P}_d}\bigg[Q(x,\rv{\xi})\bigg]
  \end{align}
where the second term models the conditional expectation of second-stage costs given a decision $x$.

\begin{example}
The rather abstract configuration introduced allows us to model various situations where the uncertain parameters (e.g., transportation costs or demands) depend on the open facilities.
Assume, for example, that $\mtc{Z}$ is a set of geographical neighborhoods.
The distribution of a customer's demand may depend, for instance, on whether there is any facility open in their neighborhood or on the distance to the nearest neighborhood with an open facility.
An example is provided in \Cref{fig:example_1}. It depicts a fictitious problem with four potential locations and one customer.
The set of locations is partitioned into three neighborhoods, $\mtc{I}_1=\{i_1,i_2\}$, $\mtc{I}_2=\{i_3\}$, $\mtc{I}_3=\{i_4\}$, thus $\mtc{Z}=\{1,2,3\}$.
This gives rise to $2^{|\mtc{Z}|}=2^3$ distributions. These are listed in \Cref{fig:example_1b}.
In the example, we see for instance that subset $\mtc{I}_1$ is activated ($\mathbb{1}_{\mtc{I}_1}(x)=1$) whenever either/both facilities $i_1$ or/and $i_2$ are open.
The different possible ways of opening facilities give rise to eight distributions, $d_1$ through $d_8$, each determined by a specific activation of subsets of $\mtc{I}$.
For instance, $d_1$ is enforced if $\mtc{I}_1$, $\mtc{I}_2$ and $\mtc{I}_3$ are all activated, which entails that facilities $i_3$ and $i_4$ are open and at least one between $i_1$ and $i_2$ is open.
Similarly, $d_2$ is enforced if $\mtc{I}_1$ and $\mtc{I}_2$ are activated but $\mtc{I}_3$ is not.
This entails that facilities $i_3$  and at least one between $i_1$ and $i_2$ are open, while $i_4$ is closed.
In this stylized example, the distribution of the demand of customer $c_1$ may for example depend on which or how many neighborhoods have open facilities.

\begin{figure}
  \centering
  
  \begin{subfigure}[b]{0.7\textwidth}
    \centering
    \tikzset{every picture/.style={line width=0.75pt}} 

\begin{tikzpicture}[x=1pt,y=1pt,yscale=-1,xscale=1]

\draw   (130.33,130.5) .. controls (130.33,125.25) and (134.59,121) .. (139.83,121) .. controls (145.08,121) and (149.33,125.25) .. (149.33,130.5) .. controls (149.33,135.75) and (145.08,140) .. (139.83,140) .. controls (134.59,140) and (130.33,135.75) .. (130.33,130.5) -- cycle ;

\draw   (160,96.17) .. controls (160,90.92) and (164.25,86.67) .. (169.5,86.67) .. controls (174.75,86.67) and (179,90.92) .. (179,96.17) .. controls (179,101.41) and (174.75,105.67) .. (169.5,105.67) .. controls (164.25,105.67) and (160,101.41) .. (160,96.17) -- cycle ;

\draw   (231.33,125.83) .. controls (231.33,120.59) and (235.59,116.33) .. (240.83,116.33) .. controls (246.08,116.33) and (250.33,120.59) .. (250.33,125.83) .. controls (250.33,131.08) and (246.08,135.33) .. (240.83,135.33) .. controls (235.59,135.33) and (231.33,131.08) .. (231.33,125.83) -- cycle ;

\draw   (207.67,186.83) .. controls (207.67,181.59) and (211.92,177.33) .. (217.17,177.33) .. controls (222.41,177.33) and (226.67,181.59) .. (226.67,186.83) .. controls (226.67,192.08) and (222.41,196.33) .. (217.17,196.33) .. controls (211.92,196.33) and (207.67,192.08) .. (207.67,186.83) -- cycle ;

\draw   (120.33,180) -- (140.67,180) -- (140.67,200.67) -- (120.33,200.67) -- cycle ;
\draw  [dash pattern={on 4.5pt off 4.5pt}] (120.71,147.01) .. controls (111.07,136.98) and (118.69,113.99) .. (137.75,95.66) .. controls (156.8,77.32) and (180.06,70.58) .. (189.71,80.61) .. controls (199.36,90.63) and (191.73,113.62) .. (172.68,131.96) .. controls (153.62,150.29) and (130.36,157.03) .. (120.71,147.01) -- cycle ;
\draw  [dash pattern={on 4.5pt off 4.5pt}] (202.16,127.6) .. controls (201.86,117.35) and (218.02,108.56) .. (238.26,107.96) .. controls (258.5,107.36) and (275.15,115.19) .. (275.45,125.44) .. controls (275.75,135.69) and (259.59,144.48) .. (239.35,145.08) .. controls (219.11,145.68) and (202.46,137.85) .. (202.16,127.6) -- cycle ;
\draw  [dash pattern={on 4.5pt off 4.5pt}] (183.42,185.2) .. controls (185.21,172.2) and (201.67,163.74) .. (220.18,166.29) .. controls (238.69,168.84) and (252.24,181.44) .. (250.45,194.44) .. controls (248.65,207.43) and (232.2,215.9) .. (213.69,213.35) .. controls (195.18,210.79) and (181.63,198.19) .. (183.42,185.2) -- cycle ;
\draw [color={rgb, 255:red, 208; green, 2; blue, 27 }  ,draw opacity=1 ] [dash pattern={on 0.84pt off 2.51pt}]  (130.67,180.33) -- (131.67,153.33) ;
\draw [color={rgb, 255:red, 208; green, 2; blue, 27 }  ,draw opacity=1 ] [dash pattern={on 0.84pt off 2.51pt}]  (140.67,180) -- (202.67,133) ;
\draw [color={rgb, 255:red, 208; green, 2; blue, 27 }  ,draw opacity=1 ] [dash pattern={on 0.84pt off 2.51pt}]  (141,191.33) -- (183.33,190.67) ;

\draw (135.33,122) node [anchor=north west][inner sep=0.75pt]  [font=\footnotesize] [align=left] {$\displaystyle i_{1}$};
\draw (165,87.67) node [anchor=north west][inner sep=0.75pt]  [font=\footnotesize] [align=left] {$\displaystyle i_{2}$};
\draw (236.33,117.33) node [anchor=north west][inner sep=0.75pt]  [font=\footnotesize] [align=left] {$\displaystyle i_{3}$};
\draw (212.67,178.33) node [anchor=north west][inner sep=0.75pt]  [font=\footnotesize] [align=left] {$\displaystyle i_{4}$};
\draw (125.67,181) node [anchor=north west][inner sep=0.75pt]  [font=\footnotesize] [align=left] {$\displaystyle c_{1}$};
\draw (111.67,84) node [anchor=north west][inner sep=0.75pt]  [font=\footnotesize] [align=left] {$\displaystyle \mathcal{I}_{1}$};
\draw (252,89.33) node [anchor=north west][inner sep=0.75pt]  [font=\footnotesize] [align=left] {$\displaystyle \mathcal{I}_{2}$};
\draw (256,168.67) node [anchor=north west][inner sep=0.75pt]  [font=\footnotesize] [align=left] {$\displaystyle \mathcal{I}_{3}$};

\end{tikzpicture}
    \caption{Locations and their partition.}
    \label{fig:example_1a}
  \end{subfigure}
  \hfill
  \begin{subfigure}[b]{0.29\textwidth}
    \centering
    \begin{tabular}{ccccc}
      \toprule
      $\mathbb{1}_{\mtc{I}_1}$&$\mathbb{1}_{\mtc{I}_2}$&$\mathbb{1}_{\mtc{I}_3}$&$d$&$\mtc{Z}_d$\\
      \midrule
      1&1&1&$d_1$&$\{1,2,3\}$\\
      1&1&0&$d_2$&$\{1,2\}$\\
      1&0&1&$d_3$&$\{1,3\}$\\
      1&0&0&$d_4$&$\{1\}$\\
      0&1&1&$d_5$&$\{2,3\}$\\
      0&1&0&$d_6$&$\{2\}$\\
      0&0&1&$d_7$&$\{3\}$\\
      0&0&0&$d_8$&$\emptyset$\\
      \bottomrule
    \end{tabular}
    \caption{Notation.}
    \label{fig:example_1b}
  \end{subfigure}

  \caption{Example partition of the set of locations and corresponding notation.}
  \label{fig:example_1}
\end{figure}
\end{example}

 The relationships between the active subsets in $\mtc{Z}$ and conditional distributions $\mathsf{P}_d$, $d\in\mtc{D}$, for $\rv{\xi}$ may be non-trivial to capture.
Several methods may apply, depending on the specific case.
For example, if similar location decisions have been made in the past, cross-sectional or panel data may be used to estimate how the relevant uncertain parameters (e.g., costs or demands) responded to those decisions across different regions. 
When demand is uncertain, customer behavior may be modeled using discrete choice models.
These estimate the probability of a (group of) customer(s) selecting a facility based on its (their) location and attributes.
Furthermore, in data-scarce environments, expert judgment may be employed to approximate how parameters might shift with the opening or closing of facilities.
This can be formalized through Bayesian methods, where expert beliefs are encoded as prior distributions, or through structured scenario analysis.
In such cases, experts provide estimates for key quantiles or ranges of the uncertain parameters under different facility configurations.
From this information probability distributions can be fitted and incorporated into the model.

\section{Solution Method}\label{sec:method}

The indicator functions $\mathbb{1}_d(x)$ in problem \eqref{eq:flp:compact} can be implemented via binary variables $\delta_d$ that take value $1$ iff distribution $d$ is enforced  and linear constraints.
Let also binary variable $y_z$ model the activation of subset $z\in\mtc{Z}$, and $\mtc{Z}_d^C$ be the complement of $\mtc{Z}_d$ to $\mtc{Z}$.
We obtain the following extensive formulation.
\begin{subequations}
  \label{eq:flp:ext}
  \begin{align}
    \label{eq:flp:ext:obj}\min &\sum_{i\in\mtc{I}}F_i x_i+\sum_{d\in\mtc{D}}\delta_d\mathbb{E}_{\mathbb{P}_d}\bigg[Q(x,\rv{\xi})\bigg]\\
    s.t.~&\label{eq:flp:ext:x-y}\sum_{i\in\mtc{I}_z}x_i\leq |\mtc{I}_z| y_z  &\forall z\in\mtc{Z}\\
                               &\label{eq:flp:ext:y-x}\sum_{i\in\mtc{I}_z}x_i\geq y_z  &\forall z\in\mtc{Z}\\
                               &\label{eq:flp:ext:d-y}\sum_{z\in\mtc{Z}_d}y_z-\sum_{z\in\mtc{Z}_d^C}y_z \geq |\mtc{Z}_d| \delta_d- |\mtc{Z}_d^C|(1-\delta_d)& \forall d \in\mtc{D}\\
                               &\label{eq:flp:ext:y-d}\sum_{z\in\mtc{Z}_d}y_z-\sum_{z\in\mtc{Z}_d^C}y_z \leq \delta_d + |\mtc{Z}_d|-1& \forall d \in\mtc{D}\\
    &x_i\in\{0,1\} & \forall i\in\mtc{I}\\
                               &y_z\in\{0,1\}&\forall z\in\mtc{Z}\\
                               &\delta_d\in\{0,1\}&\forall d \in\mtc{D}
  \end{align}
\end{subequations}
Constraints \eqref{eq:flp:ext:x-y} and \eqref{eq:flp:ext:y-x} model the implication
$$\sum_{i\in\mtc{I}_z}x_i\geq 1 \iff y_z=1$$
Constraints \eqref{eq:flp:ext:d-y} and \eqref{eq:flp:ext:y-d} model the implication
$$\sum_{z\in\mtc{Z}_d}y_z-\sum_{z\in\mtc{Z}_d^C}y_z =|\mtc{Z}_d|\iff \delta_d=1$$
That is, distribution $d$ applies if and only if the corresponding subsets of facilities ($\mtc{Z}_d$), and only those, are active.
When, for a given realization $\xi$, the second-stage problem $Q(x,\xi)$ is a (MI)LP problem and $\rv{\xi}$ is described by discrete distributions, the nonlinear term in the objective function can be linearized using standard techniques.
Hence, in this case \eqref{eq:flp:ext} can be expressed as an MILP. This allows us to use available MILP solvers to find numerical solutions.
In \Cref{sec:Appendix1} we provide an example of such a linearization for the application described in \Cref{sec:exp}.
Note, however, that the size of \eqref{eq:flp:ext} is exponential in the number of subsets $|\mtc{Z}|$.

To address the size growth, we introduce an exact decomposition method.
The method instantiates a recent extension of the classical L-Shaped method of \cite{VanW69} proposed by \cite{PanH25}.
 A brief introduction to both the methods of \cite{VanW69} and \cite{PanH25} is provided in \Cref{app:ls} to facilitate understanding for readers unfamiliar with the methodology.

The method is based on defining \textit{distribution-specific} optimality cuts.
We derive and prove validity of such cuts.
The following assumptions are made.
\begin{itemize}
\item[\textbf{A1}] For each $x$ and realization $\xi$, $Q(x,\xi)$ forms an LP.
  That is $Q(x,\xi)=\min_{w\geq 0}\{q^\top w|Ww=h-Tx\}$, where $\xi$ collects the realization $q, W, T, h$ of the random data $\rv{q},\rv{W},\rv{T},\rv{h}$. These are in turn real matrices and vectors of conformable finite dimensions.
  Hence, decisions $w$ represent second-stage operating decisions made after location decision $x$ have been made and uncertain parameters $\xi$ have materialized. 
\item[\textbf{A2}] The uncertainty $\rv{\xi}$ has discrete distributions. For each distribution $d$ there exists a finite set of scenarios $\mtc{S}_d$. Scenario $s\in\mtc{S}_d$ has probability $\pi_{sd}$ and realization $\xi_{sd}$ which collects the elements of $q_{sd},W_{sd},T_{sd},h_{sd}$.
\item[\textbf{A3}] $Q(x,\xi)$ admits feasible solution for all $x$ and $\xi$. That is, a second-stage operating decision is possible for every location decision and every realization of the parameters. As a consequence, feasibility cuts are not needed. 
\item[\textbf{A4}] There exists $\infty>U\geq \max_{d\in\mtc{D}}\max_{x\in\{0,1\}^{\vert\mtc{I}\vert}}\mathbb{E}_{\mathbb{P}_d}\bigg[Q(x,\rv{\xi})\bigg]-\min_{d\in\mtc{D}}\min_{x\in\{0,1\}^{\vert\mtc{I}\vert}}\mathbb{E}_{\mathbb{P}_d}\bigg[Q(x,\rv{\xi})\bigg]$. 
\end{itemize}
Assumptions A1, A2 and A3 are made primarily for the sake of simplicity and in coherence with the application proposed in \Cref{sec:exp}.
Assumption A4 is, instead, necessary for the proposed optimality cuts.
Essentially, it requires that second-stage costs are bounded from above and from below.
Relaxation of some of these assumptions is possible, see \cite{PanH25}.
Finally, we stress that we use the notation $Q(x,\xi)$ to refer to the second-stage problem and to its objective value for a given specification of $x$ and $\xi$ interchangeably. The meaning will be clarified when not obvious from the context.

We start by decomposing problem \eqref{eq:flp:ext} stage-wise in subproblems $Q(x,\xi)$ that address second-stage decisions,  and the following \textit{Relaxed Master Problem} (RMP) that addresses first-stage decisions.
\begin{subequations}
  \label{eq:mp}
  \begin{align}
    \label{eq:mp:obj}\min &\sum_{i\in\mtc{I}}F_i x_i+\mu\\
    s.t.~&\label{eq:mp:x-y}\sum_{i\in\mtc{I}_z}x_i\leq |\mtc{I}_z| y_z  &\forall z\in\mtc{Z}\\
                          &\label{eq:mp:y-x}\sum_{i\in\mtc{I}_z}x_i\geq y_z  &\forall z\in\mtc{Z}\\
                          &y_z\in\{0,1\}&\forall z\in\mtc{Z}
  \end{align}
\end{subequations}
Problem \eqref{eq:mp} is a relaxation of problem \eqref{eq:flp:ext}. 
A continuous decision variable $\mu$ has been introduced to replace and approximate from below the expected second-stage cost for the distributions enforced by the $x$ decisions.

Observe that the size of RMP is $\mathcal{O}(|\mtc{Z}|)$ compared to that of \eqref{eq:flp:ext} which is $\mathcal{O}(2^{|\mtc{Z}|})$. Observe further that \eqref{eq:mp} does not require indicator variables $\delta_d$.

\begin{algorithm}
\caption{Decision-Dependent L-Shaped Method for \eqref{eq:flp:ext}}\label{alg:ls}
\begin{algorithmic}[1]
\State \texttt{solved}$\gets$\texttt{false}.
\State $i\gets 0$.
\While{not \texttt{solved}}
\State Solve \eqref{eq:mp}. Let $(x^i,y^i,\mu^i)$ be its optimal solution.
\State Identify $d^i\in\mtc{D}$ such that $\mathbb{1}_{d^i}(x^i)=1$.
\State Compute $\mathbb{E}_{\mathsf{P}_{d^i}}[Q(x^i,\rv{\xi})]$.
\If{$\mu^i<\mathbb{E}_{\mathsf{P}_{d^i}}[Q(x^i,\rv{\xi})]$ }
    \State Add an optimality cut to \eqref{eq:mp} such that $\mathbb{1}_{d^i}(x)=1\implies\mu\geq o(x,y)$ and $\mu^i< o(x^i,y^i)$.
\Else   
  \State \texttt{solved}$\gets$\texttt{true}.
    \EndIf
    \State $i\gets i+1$.
\EndWhile
\end{algorithmic}
\end{algorithm}

The solution method is sketched in \Cref{alg:ls}. At iteration $i$ the algorithm solves RMP, that is, a relaxation of the original problem \eqref{eq:flp:ext}. 
Given a solution $(x^i,y^i)$ to RMP, it identifies the distribution $d^i$ enforced by $x^i$.
Following,  for decision $x^i$, it computes the expected second-stage cost under distribution $d^i$.
This is done by solving subproblem $Q(x^i,\xi_{sd^i})$ for each scenario $s\in \mtc{S}_{d^i}$ of the enforced distribution $d^i$ (and only for that distribution).
This, in turn, enables verifying whether the current value of the lower approximation, $\mu^i$, matches the second-stage expected cost conditional on $x^i$. 
If $\mu^i$ carries the correct value of the conditional expectation, the problem is solved.
Otherwise, we add a distribution-specific optimality cut $\mathbb{1}_{d^i}(x)=1\implies\mu\geq o(x,y)$ where $o(\cdot,\cdot)$ is an affine function of $x$ and $y$.
 The specification of optimality cuts is given later in \Cref{prop:oc}.
This cut makes solution $(x^i,y^i,\mu^i)$ infeasible in RMP. Essentially, it provides RMP information about the correct expected cost (under the corresponding distribution $d^i$) $\mu$ should match for solution $(x^i,y^i)$. The cut is distribution-specific in that it only applies to solutions $x$ which enforce the focal distribution $d^i$.

Observe that, since \eqref{eq:mp} is an MILP, \Cref{alg:ls} can be embedded in Branch \& Cut procedure as described by \cite{LapL93}.
Similarly, \Cref{alg:ls} is amenable to any modern acceleration technique for Benders decomposition algorithms (see e.g., the survey \citet{RahCGR17}).
We present this basic version for the sake of simplicity.

We proceed by describing how to obtain optimality cuts.
 To that effect, observe that a direct consequence of A1, A3 and A4 is that $Q(x,\xi)$ is a feasible linear program with a bounded objective value.
As such, it has a strong dual, say $Q^D(x,\xi)$ with $Q^D(x,\xi)=Q(x,\xi)$ for all $x$ and $\xi$. Dual solutions are required for building optimality cuts, as we see next. 
\begin{proposition}\label{prop:oc}
  Let $(x^i,y^i,\mu^i)$ a solution to RMP. Let $d^i$ be such that $\mathbb{1}_{d^i}(x^i)=1$.
  Assume $\mu^i<\mathbb{E}_{\mathsf{P}_{d^i}}\left[Q(x^i,\rv{\xi})\right]$.
  Then, $(x^i,y^i,\mu^i)$ violates the following inequality.
  \begin{equation}\label{eq:oc}
    \mu\geq \sum_{s\in\mtc{S}_{d^i}}\pi_{sd}\left[(\rho^i_{sd})^\top\left(h_{sd}-T_{sd}x\right)\right]-U\bigg(|\mtc{Z}_{d^i}|-\sum_{z\in\mtc{Z}_{d^i}}y_z+\sum_{z\in\mtc{Z}^C_{d^i}}y_z\bigg)
  \end{equation}
  Where $\rho_{sd}^i$ is an optimal solution to the dual of $Q(x^i,\xi_{sd})$.

\end{proposition}
  \begin{proof}{Proof}
    By assumption $\mu^i<\mathbb{E}_{\mathsf{P}_{d^i}}\left[Q(x^i,\rv{\xi})\right]$. Therefore, we can write
    \begin{align*}
      \mu^i<&\mathbb{E}_{\mathsf{P}_{d^i}}\left[Q(x^i,\rv{\xi})\right]\\
            =&\sum_{s\in\mtc{S}_{d^i}}\pi_{sd^i}\left[Q(x^i,\xi_{sd^i})\right]\\
            =&\sum_{s\in\mtc{S}_{d^i}}\pi_{sd^i}\left[(\rho^i_{sd^i})^\top\left(h_{sd^i}-T_{sd^i}x^i\right)\right]\\
            =&\sum_{s\in\mtc{S}_{d^i}}\pi_{sd^i}\left[(\rho^i_{sd^i})^\top\left(h_{sd^i}-T_{sd^i}x^i\right)\right]-U\bigg(|\mtc{Z}_{d^i}|-\sum_{z\in\mtc{Z}_{d^i}}y^i_z+\sum_{z\in\mtc{Z}^C_{d^i}}y^i_z\bigg)
    \end{align*}
    The first equality follows from assumption A2. The second equality holds by assumption A3 and strong duality of $Q(x^i,\xi_{sd})$. The third equality holds because $d^i$ implies $\sum_{z\in\mtc{Z}_{d^i}}y^i_z=|\mtc{Z}_{d^i}|$ and $\sum_{z\in\mtc{Z}^C_{d^i}}y^i_z=0$ and proves the inequality violated as required. 
  \end{proof}

  The next proposition clarifies that cuts \eqref{eq:oc} are not violated by solutions to RMP where $\mu$ bears the correct value of the second-stage expectation.
  
\begin{proposition}\label{prop:oc:safe}
  Let $(x^i,y^i,\mu^i)$ a solution to RMP. Let $d^i$ be such that $\mathbb{1}_{d^i}(x^i)=1$.
  Assume $\mu^i=\mathbb{E}_{\mathsf{P}_{d^i}}\left[Q(x^i,\rv{\xi})\right]$.
  Then, $(x^i,y^i,\mu^i)$ satisfies inequality \eqref{eq:oc} generated for some $d^v$ and $\rho^v=(\rho_{sd^v}^v)_{s\in\mtc{S}_{d^v}}$.

\end{proposition}

  \begin{proof}{Proof}
    We first prove the statement for the case when $d^i=d^v$ and then for the case when $d^i\neq d^v$.
    
    Assume $d^i=d^v$. Then
    \begin{align*}
      \mu^i=&\mathbb{E}_{\mathsf{P}_{d^i}}\left[Q(x^i,\rv{\xi})\right]\\
            =&\sum_{s\in\mtc{S}_{d^i}}\pi_{sd^i}\left[Q(x^i,\xi_{sd^i})\right]\\
            =&\sum_{s\in\mtc{S}_{d^i}}\pi_{sd^i}\left[(\rho^i_{sd^i})^\top\left(h_{sd^i}-T_{sd^i}x^i\right)\right]\\
            \geq&\sum_{s\in\mtc{S}_{d^i}}\pi_{sd}\left[(\rho^v_{sd^i})^\top\left(h_{sd^i}-T_{sd^i}x^i\right)\right]\\
            =&\sum_{s\in\mtc{S}_{d^i}}\pi_{sd^i}\left[(\rho^v_{sd^i})^\top\left(h_{sd^i}-T_{sd^i}x^i\right)\right]-U\bigg(|\mtc{Z}_{d^i}|-\sum_{z\in\mtc{Z}_{d^i}}y^i_z+\sum_{z\in\mtc{Z}^C_{d^i}}y^i_z\bigg)
    \end{align*}
    The central inequality follows from dual optimality of $\rho^i_{sd^i}$ and dual feasibility of $\rho^v_{sd^i}$.
    The last equality holds because $d^i$ implies $\sum_{z\in\mtc{Z}_{d^i}}y^i_z=|\mtc{Z}_{d^i}|$ and $\sum_{z\in\mtc{Z}^C_{d^i}}y^i_z=0$ and proves the inequality satisfied as required.

    Assume $d^i\neq d^v$. Then
    \begin{align*}
      \mu^i= \mathbb{E}_{\mathsf{P}_{d^i}}\left[Q(x^i,\rv{\xi})\right]\geq&\mathbb{E}_{\mathsf{P}_{d^v}}\left[Q(x^i,\rv{\xi})\right]-U = \sum_{s\in\mtc{S}_{d^v}}\pi_{sd^v}\left[(\rho^v_{sd^v})^\top\left(h_{sd^v}-T_{sd^v}x^i\right)\right]-U\\
      \geq &\sum_{s\in\mtc{S}_{d^v}}\pi_{sd^v}\left[(\rho^v_{sd^v})^\top\left(h_{sd^v}-T_{sd^v}x^i\right)\right]-U\bigg(|\mtc{Z}_{d^v}|-\sum_{z\in\mtc{Z}_{d^v}}y^i_z+\sum_{z\in\mtc{Z}^C_{d^v}}y^i_z\bigg)
    \end{align*}
    The first inequality follows from the definition of $U$. The second equality holds by strong duality.
    Finally, the last inequality holds because $|\mtc{Z}_{d^v}|-\sum_{z\in\mtc{Z}_{d^v}}y^i_z+\sum_{z\in\mtc{Z}^C_{d^v}}y^i_z\geq 1$ and proves the inequality satisfied as required.
  \end{proof}

\begin{proposition}\label{prop:conv}
  \Cref{alg:ls} converges in a finite number of iterations.

\end{proposition}
  \begin{proof}
    Finite convergence of the algorithm follows immediately from the finite number of optimality cuts that it is possible to generate.
    This, in turn, follows from the finite number of distributions in $\mtc{D}$ and the finite number of extreme points $\rho_{sd}$ in the dual to the problems $Q(\cdot,\xi_{sd})$ (observe that the feasible region of the dual to $Q(\cdot,\xi_{sd})$ does not depend on the decision variables of RMP since they appear in the right-hand side of the constraints of $Q(x,\xi)$, see Assumption A1). When, at most, all the finitely many cuts are addedd to RMP, RMP necessarily produces the optimal solution, see \citet{PanH25}.
  \end{proof}

  In summary, \Cref{prop:oc} clarifies that optimality cuts render infeasible in RMP any solution $(x,y,\mu)$ where $\mu$ is a strict lower-approximation of the expected second-stage cost.
  \Cref{prop:oc:safe} clarifies that the cuts would instead preserve solutions $(x,y,\mu)$ to RMP, where $\mu$ bears the correct value of the expected second-stage cost.
  Among these solutions we find the optimal solution.
  Finally, \Cref{prop:conv} clarifies that the method will stop generating cuts and, hence, deliver the optimal solution in finitely many iterations. 
  
\section{Numerical Study and Problem-Specific Enhancements}\label{sec:exp}

We present the results of experiments on a facility location problem with endogenous stochastic demand.
The problem is adapted from \cite{LapL94} and consists of optimally determining the location and size of facilities to serve uncertain future demand.
Since demand is uncertain, the requirement that demands should be met in all circumstances becomes unrealistic \cite{Lou86,LapL94}.
For this reason, a selling price is introduced. The problem becomes that of finding a trade-off between initial investments and future profits.
Opening too few facilities may result in lost sales. Opening too many facilities may result in unused resources.
The problem is further complicated by the influence of location decisions on the distribution of demands.
This relaxes the classical assumption that demand is independent of distances, see e.g., \cite{Lou86}.

The notation introduced above is augmented as follows.
We let $C_i$ be the maximum capacity of a facility located at $i\in\mtc{I}$.
We let $\mtc{J}$ be the set of customer locations.
For all $i\in\mtc{I}$ and $j\in\mtc{J}$, we let $R_{ij}$ be the net revenue per unit of demand of customer $j$ satisfied by facility $i$.
Given an underlying measurable space $(\Omega,\mtc{F})$, we let $\rv{\xi}_j:\Omega\to\mathbb{R}_+$ be a random variable representing the uncertain demand of customer $j\in\mtc{J}$.
Contrary to standard assumptions, the distribution of $\rv{\xi}=(\rv{\xi}_j)_{j\in\mtc{J}}$ depends on location decisions as we see next.

We assume the probability distribution of customer demands $\rv{\xi}_j$ is fully determined by the open facilities.

To mimic a phenomenon by which demand is influenced by proximity to, and/or density of, open facilities we proceed as follows.
We give set $\mtc{Z}$ the interpretation of a set of geographical zones.
We define $z^{(n)}(j)$ to be the $n$-th zone closest to customer $j$, with $n=1,\ldots, |\mtc{Z}|$.
We let the demand of customer $j$ follow a normal distribution truncated to the left at $0$ and parametrized by indicator variables $\mathbb{1}_{\mtc{Z}_d}(z)$.
Particularly, we let the mean demand for customer $j$ under distribution $d$ be defined as
%
$$\mu_{jd}=g(\mathbb{1}_{\mtc{Z}_d}(1),\cdots,\mathbb{1}_{\mtc{Z}_d}(|\mtc{Z}|)):=\mu_j+\sum_{n=1}^{|\mtc{Z}|}(-1)^{s(n)}\alpha_{nj}\mu_j\mathbb{1}_{\mtc{Z}_d}(z^{(n)}(j))$$
with $\alpha_{nj}\in[0,1)$ and $s(n):\mathbb{Z}\to \{0,1\}$. That is, the mean under distribution $d$ is equal to some base-level mean $\mu_j$,
which increases or decreases (depending on $s(n)$) by $100\alpha_{nj}\%$ whenever there is an open facility in the $n$-th furthest zone.
In a similar way, the standard deviation is defined as 
$$\sigma_{jd}=h(\mathbb{1}_{\mtc{Z}_d}(1),\cdots,\mathbb{1}_{\mtc{Z}_d}(|\mtc{Z}|)):=\sigma_j-\sum_{n=1}^{|\mtc{Z}|}(-1)^{s(n)}\beta_{nj}\sigma_j\mathbb{1}_{\mtc{Z}_d}(z^{(n)}(j))$$
with $\sigma_j$ being a base-level standard deviation and $\beta_{nj}\in[0,1)$.
By adjusting the $\alpha_{nj}$ and $\beta_{nj}$ coefficients one can model different types of interactions between locations and demand distributions.
The instances presented in \Cref{sec:exp:inst} describe four such interactions. 
For the sake of simplicity we let customer demands be uncorrelated.
We emphasize that the choice of a linear relationship between the open zones and the parameters of the distribution, as well as the specifications provided in \Cref{sec:exp:inst}, serve as a proxy and are not intended to accurately represent real-world dynamics.

Given this notation we can formulate the \textit{Capacitated Facility Location Problem under Endogenous Uncertain Demand} as follows.
\begin{subequations}
  \label{eq:flpsd:compact}
  \begin{align}
    \label{eq:flpsd:compact:obj}\min_{x\in\{0,1\}^{|\mtc{I}|}}  &F^\top x-\sum_{d\in\mtc{D}}\mathbb{1}_d(x)\mathbb{E}_{\mathsf{P}_d}\bigg[Q(x,\rv{\xi})\bigg]
  \end{align}
  where
  \begin{align}
    \label{eq:flpsd:2st:obj}Q(x,\xi) = & \max \sum_{i\in\mtc{I}}\sum_{j\in\mtc{J}}R_{ij}w_{ij}\\
    \label{eq:flpsd:2st:1}\text{s.t.~}& \sum_{i\in\mtc{I}}w_{ij}\leq \xi_j & j\in\mtc{J}\\    
    \label{eq:flpsd:2st:2}        &\sum_{j\in\mtc{J}}w_{ij}\leq C_ix_i& i\in\mtc{I}\\
    \label{eq:flpsd:2st:3}        &w_{ij}\geq 0&i\in\mtc{I},j\in\mtc{J}
  \end{align}
\end{subequations}
The objective of the problem is that of opening facilities in order to maximize conditional expected profits.
The second-stage problems maximize revenues, with second-stage constraints ensuring that facility capacities are respected and sales do not exceed demand. 
Variables $w_{ij}$ represent the amount of demanded good from customer $j$ satisfied by facility $i$.
The extensive linearized model for this problem is provided in \Cref{sec:Appendix1}.

\subsection{A Problem-Specific Valid Inequality}\label{sec:exp:vi}

By careful analysis of $Q(x,\xi)$ it is possible to obtain a valid inequality.

\begin{proposition}\label{prop:vi}
  
  \begin{equation}\label{eq:vi}
    \mu\leq \sum_{i\in\mtc{I}}\sum_{j\in\mtc{J}}R_{ij}\min\{C_i,\max_{d\in\mtc{D}}\mathbb{E}_{\mathsf{P}_d}[\rv{\xi}_j]\}x_i
  \end{equation}
  is a valid inequality for the RMP of \eqref{eq:flpsd:compact}.
\end{proposition}
\begin{proof}
  We need to show that, for all $x^v\in\{0,1\}^{|\mtc{I}|}$, we have
  $$\mathbb{E}_{\mathsf{P}_{d^v}}\bigg[Q(x^v,\rv{\xi})\bigg]\leq \sum_{i\in\mtc{I}}\sum_{j\in\mtc{J}}R_{ij}\min\{C_i,\max_{d\in\mtc{D}}\mathbb{E}_{\mathsf{P}_d}[\rv{\xi}_j]\}x_i^v$$
  Take $x^v\in\{0,1\}^{|\mtc{I}|}$ and let $d^v$ be the enforced distribution (i.e., $\mathbb{1}_{d^v}(x^v)=1$).
  We start by observing that $Q(x^v,\xi)$ is concave in $\xi$. Hence, by Jensen's inequality we can state that
  $$\mathbb{E}_{\mathsf{P}_{d^v}}\bigg[Q(x^v,\rv{\xi})\bigg]\leq Q(x^v,\mathbb{E}_{\mathsf{P}_{d^v}}[\rv{\xi}])$$
  Let $w^*$ be an optimal solution to $Q(x^v,\mathbb{E}_{\mathsf{P}_{d^v}}[\rv{\xi}])$.
  Given constraints \eqref{eq:flpsd:2st:1} and \eqref{eq:flpsd:2st:2}, it is easy to see that, in \eqref{eq:flpsd:2st:obj}, for each $i\in\mtc{I}$ and $j\in\mtc{J}$ we have
  $$R_{ij}w^*_{ij}\leq R_{ij}\min\{C_i,\mathbb{E}_{\mathsf{P}_{d^v}}[\rv{\xi}_j]\}x_i^v$$
  That is, when $x^v_i=0$ the flow between $i$ and $j$ is $0$, while when $x_i^v=1$ the flow is bounded above by the smallest between the capacity of $i$ and the demand of $j$. Following
  $$Q(x^v,\mathbb{E}_{\mathsf{P}_{d^v}}[\rv{\xi}])\leq \sum_{i\in\mtc{I}}\sum_{j\in\mtc{J}}R_{ij}\min\{C_i,\mathbb{E}_{\mathsf{P}_{d^v}}[\rv{\xi}_j]\}x_i^v\leq \sum_{i\in\mtc{I}}\sum_{j\in\mtc{J}}R_{ij}\min\{C_i,\max_{d\in\mtc{D}}\mathbb{E}_{\mathsf{P}_d}[\rv{\xi}_j]\}x_i^v$$
  This completes the proof. 
\end{proof}

Observe that the inequality is reversed by effect of the maximization in the subproblem.
Observe also that \eqref{eq:vi} is, in general, not tight. That is, unlike optimality cuts \eqref{eq:oc}, it does not guarantee that, for the optimal solution $x^*$ (with distribution $d^*$), one has 
$$\mu^*=\mathbb{E}_{\mathsf{P}_{d^*}}\bigg[Q(x^*,\rv{\xi})\bigg]$$
Hence, convergence of \Cref{alg:ls} is not guaranteed by \eqref{eq:vi} only, but by optimality cuts \eqref{eq:oc}.
On the other hand, in \Cref{sec:exp:res} we show that valid inequality \eqref{eq:vi} helps reducing the optimality gap and finding feasible solutions for large instances.

\subsection{Instances}\label{sec:exp:inst}

Instances of the problem are generated by adapting and expanding the recipe provided by \citet{LapL94}.
We start by generating $|\mtc{J}|$ customers randomly with uniform probability in a square of area $[0,100]\times [0,100]$.
Customer demands $\rv{\xi}_j$ are assumed to follow uncorrelated truncated normal distributions parametrized by location decisions.
Each customer $j\in\mtc{J}$ has base mean $\mu_j$ randomly drawn uniformly in the interval $[10,50]$ and ratio $\sigma_j/\mu_j$ randomly drawn uniformly in $[0.05,0.35]$.
Three potential facility locations are then generated in the $10\times 10$ square centered around each of the three customers with the largest $\mu_j$.
The remaining $|\mtc{I}|-3$ locations are randomly generated in the $60\times 60$ square centered around the point $(50,50)$.

To create geographical zones $\mtc{Z}$ we partition facility locations into $|\mtc{Z}|$ clusters according to their coordinates.
For each customer, we rank zones $\mtc{I}_z$ in non-decreasing order of the distance between the customer and the centroid of the zone.
Thus, we let $z^{(n)}(j)$, $n=1,\ldots,\vert\mtc{Z}\vert$ be computed from the distance between customer $j$ and the centroids of the clusters.
We test the following configurations of the $\alpha_{nj}$ and $\beta_{nj}$ values.
  \begin{itemize}
  \item[A] In the first configuration we set $s(n)=0$,  $\alpha_{nj}=0.5^n$ and $\beta_{nj}=0.4^n$, $n=1,\ldots,\vert\mtc{Z}\vert$, $j\in\mtc{J}$.
    This entails that the contribution of the the open zones to mean and standard deviation decreases exponentially with the rank of their distance from the customer.
    Hence, demand is influenced by both proximity and density of open zones. 
  \item[B] In the second configuration we set $s(n)=0$, $n=1,\ldots,\vert\mtc{Z}\vert$, $\alpha_{1j}=0.5$, $\beta_{1j}=0.4$ and $\alpha_{nj}=\beta_{nj}=0$ for $n=2,\ldots,\vert\mtc{Z}\vert$ and $j\in\mtc{J}$.
    This entails that the distribution of the demand of customer $j$ is influenced only by the opening of facilities in the zone closest to customer $j$, $z^{(1)}(j)$, and not influenced by the opening of other facilities. 
  \item[C] In the third configuration we let $\bar{n}_{jd}=\arg\min_{n=1,\ldots,|\mtc{Z}|}\{n\vert \mathbb{1}_{\mtc{Z}_d}(z^{(n)}(j))=1\}$ identify the order of the nearest open zone under distribution $d$.
    Then we set $s(n)=0$, $n=1,\ldots,\vert\mtc{Z}\vert$,  $\alpha_{nj}=0.5^n\mathbb{1}_{n=\bar{n}_{jd}}$ where $\mathbb{1}_{n=\bar{n}_{jd}}$ is the indicator of the event $n=\bar{n}_{jd}$. Hence, the mean of the demand increases by $0.5^n$ if and only if
    the closest open zone is the $n$-furthest. In a similar way $\beta_{nj}=0.4^n\mathbb{1}_{n=\bar{n}_{jd}}$. 
  \item[D] In the fourth configuration we set $s(n)=0$, for $n=1$ and $s(n)=1$ for $n=2,\ldots,\vert\mtc{Z}\vert$, and $\alpha_{nj}=0.5^n$ and $\beta_{nj}=0.4^n$ $n=1,\ldots,\vert\mtc{Z}\vert$, $j\in\mtc{J}$.
    This entails that the opening of facilities in the closest zone increases the mean demand and decreases the standard deviation.
    On the contrary, the opening of facilities in zones further away than the closest decreases the demand and increases the standard deviation. This configuration mimics a scenario where demand decreases in certain areas if facilities are located elsewhere. 
\end{itemize}
We refer to these configurations as \textit{demand types} ($\texttt{dt}$) A, B, C and D.
 We note that these configurations are merely stylized representations of real-life dynamics, which are necessarily more complex. In particular, the proposed configurations serve only as test cases for the methodology presented. 
Given a distribution, realizations of the demands $\rv{\xi}_j$ are sampled iid from the underlying truncated normal distributions.

We test instances that vary in the number of facilities, customers, zones and scenarios as well as in the parameters of the problem.
Particularly, we test instance sizes with $(|\mtc{I}|,|\mtc{J}|,|\mtc{Z}|,|\mtc{S}_d|)\in\{10,15,20,25\}\times \{50,100\}\times \{5,7,10\}\times \{50,100\}$.
For each of the resulting combinations we test seven configurations of the following parameters
\begin{itemize}
\item $R_{ij}=R$, the net revenue per unit of demand satisfied (we assume it is the same for all customers),
\item $C_i=C|\mtc{J}|$, the capacity per facility (we assume all facilities have the same capacity),
\item $F_i=O|\mtc{J}|$, the fixed opening cost per facility (we assume all facilities have the same opening cost).
\end{itemize}

The seven configurations of the values for constants $C$, $O$ and $R$ are provided in \Cref{tab:lapParams}.

Finally, since instances are randomly generated, for every possible configuration of the control parameters we generate three random instances.
This results in $48\times 7\times 3=1008$ instances for each demand type.

\begin{table}
  \centering
  \caption{Configurations of the parameters in the instances.}\label{tab:lapParams}
  \begin{tabular}{cccc}
    \toprule
    Configuration & $C$ & $O$ & $R$\\
    \midrule
    1 &  15	& 500	&400\\
    2 &  12.5	& 500	&400\\
    3 & 17.5	& 500	&400\\
    4&  15	& 250	&400\\
    5&  15	& 750	&	400\\
    6 & 15	& 500	&	200\\
    7 & 15	 &500		&600\\
    \bottomrule
  \end{tabular}
\end{table}

\subsection{Results}\label{sec:exp:res}

The algorithm has been implemented in Java using Cplex 20.1 libraries with optimality cuts added only to integer nodes in the Branch \& Bound tree (i.e., we start the Branch \& Bound method only once for RMP and separate cuts upon reaching integer nodes).
The extensive linearized form of the problem (see \Cref{sec:Appendix1}) has likewise been solved with Cplex 20.1.
All tests have been run on machines with 40 CPUs and 188 GB of memory.
No further efficiency measures were implemented in our algorithm. Particularly, parallel search was turned off due to the presence of callbacks (see \cite{Cplex20.1}) and no further parallelization was implemented.
Hence, our algorithm run entirely in a single thread.
On the other hand, when solving the extensive form, Cplex was allowed to use parallel search (with up to 32 threads) in the Branch \& Bound algorithm.
All tests have been run with Cplex's default tolerances, including the optimality tolerance of $10^{-4}$.
This is calculated as
$|\mathtt{bestbound}-\mathtt{bestinteger}|/(10^{-10}+|\mathtt{bestinteger}|)$. The time limit was set to $1800$ seconds in all tests, unless otherwise specified.
 All results presented are obtained without the use of valid inequality \eqref{eq:vi}, unless explicitly stated otherwise.

For both our algorithm (LS) and the extensive linearized form (Cplex), we report statistics on the results of the tests in \Cref{tab:byFacility}.
Particularly, we report on the following metrics:
\begin{itemize}
\item The average optimality gap (column \texttt{avgGap}) as a percentage. This is calculates as $100*|\mathtt{bestbound}-\mathtt{bestinteger}|/(10^{-10}+|\mathtt{bestinteger}|)$.
\item The number of instances for which a feasible solution was found (column \texttt{feas}).
\item The number of instances solved to a gap smaller than $0.5$\% (column $<0.5\%$).
\item The average solution time (column \texttt{avgT}) in seconds.
\item The minimum optimality gap (column \texttt{minGap}) among all feasible solutions found (for Cplex only).
\end{itemize}

\begin{table}
  \centering
  \caption{Results summarized by number of facilities and demand type (\texttt{dt}).}  \label{tab:byFacility}
\begin{tabular}{ll|ccccc|cccc}
  \toprule
&   &  \multicolumn{5}{c}{Cplex}  & \multicolumn{4}{c}{LS}    \\
  $|\mtc{I}|$ & \texttt{dt} & \texttt{minGap}[\%] & \texttt{avgGap}[\%] &  \texttt{feas}&  $<0.5\%$ &  \texttt{avgT}[sec] &   \texttt{avgGap}[\%] & \texttt{feas} & $<0.5\%$ & \texttt{avgT}[sec]  \\
    \midrule
10 &A& 0.00 & $>10^3$ & 147 & 6 & 1777.26 & 0.00 & 252 & 252 & 95.21 \\
15 &A& 525.71 & $>10^3$ & 147 & 0 & 1814.87 & 0.00 & 252 & 252 & 113.31 \\
20 &A& 992.99 & $>10^3$& 147 & 0 & 1821.34 & 0.00 & 252 & 252 & 128.21 \\
25 &A& 1296.08 & $>10^3$ & 138 & 0 & 1975.10 & 0.00 & 252 & 252 & 144.27 \\
  \midrule
10  &B & 0.00 & $>10^3$ & 168 & 14 & 1866.47 & 0.00 & 252 & 252 & 70.74 \\
15 &B & 471.84 & $>10^3$ & 146 & 0 & 1803.70 & 0.00 & 252 & 252 & 81.31 \\
20 & B& 759.11 & $>10^3$ & 147 & 0 & 1824.42 & 0.00 & 252 & 252 & 89.66 \\
25  & B& 1038.89 & $>10^3$ & 144 & 0 & 1909.24 & 0.00 & 252 & 252 & 99.85 \\
  \midrule
10 &C& 0.00 & $>10^3$ & 168 & 17 & 1851.07 & 0.00 & 252 & 252 & 72.39 \\
15 &C& 508.13 & $>10^3$ & 146 & 0 & 1803.89 & 0.00 & 252 & 252 & 81.31 \\
20 &C& 978.14 & $>10^3$ & 147 & 0 & 1824.64 & 0.00 & 252 & 252 & 90.82 \\
  25 &C& 1144.77 & $>10^3$ & 143 & 0 & 1906.63 & 0.00 & 252 & 252 & 101.66 \\
  \midrule
   10& D &  0.00 &  $>10^3$ &  167  &  21 &  1854.23 &  0.00 &  252 &  252 &  96.07 \\
   15 & D&  374.57 &  $>10^3$ &  141 &  0 &  1804.57 &  0.00 &  252 &  252 &  107.83 \\
   20 & D&  741.84 &  $>10^3$ &  142 &  0 &  1808.57 &  0.00 &  252 &  252 &  116.66 \\
   25 & D&  1024.08 &  $>10^3$ &  133 &  0 &  1901.00 &  0.00 &  252 &  252 &  128.79 \\
  \bottomrule
  \end{tabular}
\end{table}

The statistics in \Cref{tab:byFacility} illustrate that solving the extensive linearized formulation is viable only for the smallest instances (i.e., with $10$ facilities).
Also in that case, the solver solves to (near) optimality only a small fraction of the instances. 
On the remaining instances, Cplex finds feasible solutions to approximately one-half to two-thirds of the instances. However, for the majority of these it reports a three-digit optimality gap when the time limit expires.
Hence, the exponential size of the formulation makes it applicable only to problems of very limited size.

Our method, on the other hand, performs significantly better.
All tested instances were solved to optimality (within the default $10^{-4}$ tolerance).
The average solution time is also reasonably small. It never exceeds $4$ minutes. Furthermore, the solution time appears to be growing sub-linearly with the number of facilities.

\begin{table}
  \centering
  \caption{Results summarized by number of zones and demand type (\texttt{dt}).}  \label{tab:byZones}
  \begin{tabular}{ll|ccccc|cccc}
    \toprule
     &   & \multicolumn{5}{c}{Cplex}  & \multicolumn{4}{c}{LS}    \\
    $|\mtc{Z}|$ & \texttt{dt}& \texttt{minGap}[\%] & \texttt{avgGap}[\%] &  \texttt{feas}&  $<0.5\%$ &  \texttt{avgT}[sec] &   \texttt{avgGap}[\%] & \texttt{feas} & $<0.5\%$ & \texttt{avgT}[sec]  \\
    \midrule
    5       &   A  &     0.00 &                                            $>10^3$ &           333 &             6 &   1809.79 &        0.00 &         336 &         336 &   13.55 \\
    7       &   A  &   583.34 & $>10^3$ &           246 &             0 &   1893.02 &        0.00 &         336 &         336 &   43.34 \\
    10      &   A  &      - &                                                - &             0 &             0 &       - &        0.00 &         336 &         336 &  303.86 \\
    \midrule
    5 & B&0.00 &  $>10^3$ & 334 & 13 & 1781.08 & 0.00 & 336 & 336 & 9.86\\
    7 & B& 0.00 & $>10^3$ & 270 & 1 & 1938.32 & 0.00 & 336 & 336 & 30.82 \\
    10 & B& - & - & 0 & 0 & - & 0.00 & 336 & 336 & 215.48 \\
    \midrule
    5 & C&0.00 &  $>10^3$ & 335 & 17 & 1773.70 & 0.00 & 336 & 336 & 10.52 \\
    7 & C&171.12 & $>10^3$ & 269 & 0 & 1936.90 & 0.00 & 336 & 336 & 31.85 \\
    10 &C &- & - & 0 & 0 & - & 0.00 & 336 & 336 & 217.26 \\
    \midrule
    5 &D& 0.00 & $>10^3$ & 336 & 20 & 1775.32 & 0.00 & 336 & 336 & 13.55 \\
    7 &D& 0.00 & $>10^3$ & 247 & 1 & 1932.15 & 0.00 & 336 & 336 & 40.43 \\
    10&D & - & - &  0 &  0 & - & 0.00 & 336 & 336 & 283.03 \\
\bottomrule
  \end{tabular}
\end{table}

Further evidence on the performance of the method is provided in \Cref{tab:byZones}. Here the results are aggregated by the number of zones.
Recall that the number of distributions is exponential in the number of zones.
It can be noted that solving the extensive formulation is viable only for the smallest number of zones considered ($|\mtc{Z}|=5$).
In this case, the solver found a feasible solution to almost all the $336$ instances of that size.
Nevertheless, only a small fraction of these instances were solved to a gap smaller than $0.5\%$.
The average optimality gap within the time limit remains larger than $10^3$\%.

\color{black}
Observe, in \Cref{tab:byFacility} and \Cref{tab:byZones} that average solution times (column ``\texttt{avgT}'') for ``Cplex'' were calculated by taking the average solution time of both the instances that were solved to optimality and those that were not solved to optimality. For the latter the solver terminated when reaching the time limit of 1800 seconds, or a slightly longer times dure to system overhead. This implies that the solution time of the instanced that were solved to optimality by Cplex is potentially slightly overestimated. The average solution time of all instances that where solved to optimality by the solver is 1242.78 seconds.
\color{black}

As seen in \Cref{tab:byFacility}, our method solves all the instances.
\Cref{tab:byZones} offers some more insights on the solution time.
It can be noticed that the solution time grows significantly with the number of zones.
This behavior is expected due to the exponential growth of the number of distributions.
This entails that the method needs to generate an amount of optimality cuts sufficient to approximate exponentially more distributions.
However, we observed that the average number of optimality cuts per distribution is consistently close to $1$, the mode is always $1$, and the standard deviation smaller than $0.20$.
This suggests that the method is usually able to obtain enough information on each sub-region of the partition of $\mtc{I}$ with only one optimality cut, and typically not many more than that.
Furthermore, the average number of Branch \& Bound nodes explored by our method is below two thousand nodes even for the largest instances, see \Cref{tab:nodes}.
This translates into a relatively low memory requirement of the Branch \& Bound tree, which never exceeded $0.1MB$.
Consequently, even in the case with $2^{|\mtc{Z}|}=2^{10}$ distributions, the solution time is, on average, approximately five minutes. For the same instances, the solver could not find feasible solutions.

\begin{table}
  \centering
  \caption{Average number of Branch \& Bound nodes explored aggregated by number of zones and demand type (\texttt{dt}).}  \label{tab:nodes}
    
\begin{tabular}{ccc}
\toprule
$|\mtc{Z}|$& \texttt{dt} & Average \# nodes explored                \\
\midrule
5       & A &         68.42 \\
7       &  A &      273.07 \\
  10      & A &      1937.21 \\
  \midrule
   5 &B& 73.84 \\
 7 &B& 253.56 \\
 10& B& 1945.53 \\
  \midrule
 5 & C&128.78 \\
7 & C& 273.89 \\
  10 &C& 1966.99 \\
  \midrule
  5 &D& 150.82 \\
7 &D& 285.43 \\
  10&D & 1958.94 \\
\bottomrule
\end{tabular}
\end{table}

As we have seen in \Cref{tab:byZones}, the number of zones is the input that impacts the solution time the most.
To assess the limits of the proposed method, we performed experiments with a higher number of zones.
Particularly, for simplicity, we focused only on the instances with $|\mtc{I}|=20$, $|\mtc{J}|=50$, $|\mtc{S}|=50$.
We then increased the number of zones up to and including the limiting case $|\mtc{Z}|=|\mtc{I}|$ (i.e., where every location decision determines a different distribution).
This setup generates instances with up to $2^{20}$ distributions. 
The time limit is set to $36 000$ seconds ($10$ hours) in preparation for longer computation times. 
\Cref{tab:byZonesLarge} reports the results with and without valid inequality \eqref{eq:vi}.

\begin{table}
  \centering
  \caption{Results summarized by number of zones for all instances with $|\mtc{I}|=20$, $|\mtc{J}|=50$, $|\mtc{S}|=50$, with and without valid inequality. The time limit is set to $36 000$ seconds.}  \label{tab:byZonesLarge}
    
\begin{tabular}{l|rrrr|rrrr}
\toprule
 & \multicolumn{4}{c}{LS} & \multicolumn{4}{c}{LS+VI} \\
$|\mtc{Z}|$ & \texttt{avgGap}[\%] & \texttt{feas} & $<0.5\%$ & \texttt{avgT}[sec] &  \texttt{avgGap}[\%] & \texttt{feas} & $<0.5\%$ & \texttt{avgT}[sec]\\
  \midrule
  12 & 0.00 & 84 & 84 & 431.00 & 0.00 & 84 & 84 & 440.21 \\
15 & 0.00 & 84 & 84 & 4257.74 & 0.00 & 84 & 84 & 4429.29 \\
17 & 293.45 & 82 & 74 & 27138.76 & 264.18 & 83 & 77 & 26175.02 \\
20 & 7954.73 & 61 & 0 & 36261.87 & 7598.54 & 71 & 0 & 36148.41 \\
\bottomrule
\end{tabular}
\end{table}

The proposed method was able to solve to optimality all instances with up to $15$ zones (i.e., $2^{15}$ distributions) within less than two hours.
When the number of zones is higher the quality of the results inevitably worsens.
Nevertheless, with $|\mtc{Z}|=17$ the method still finds a feasible solution to $82$ out of the $84$ instances of that size.
For $74$ of these instances, the solution has an optimality gap smaller than $0.5\%$.
With $|\mtc{Z}|=20$ the method was not able to solve any instance to optimality within the given time limit.
However, it found a feasible solution for $61$ out of the $84$ instances of that size, though with a high optimality gap.
The addition of valid inequality \eqref{eq:vi} has an impact on the largest instances.
When $|\mtc{Z}|=17$ it allows the method to terminate with a feasible solution to three additional instances and to close the gap of one additional instance.
When $|\mtc{Z}|=20$ the valid inequality allows finding more feasible solutions, though the optimality, even if slightly reduced, remains high.

\section{Potential Impact and Sensitivity Analysis}\label{sec:ex}
\begin{figure}
  \centering
  \includegraphics[width=\textwidth]{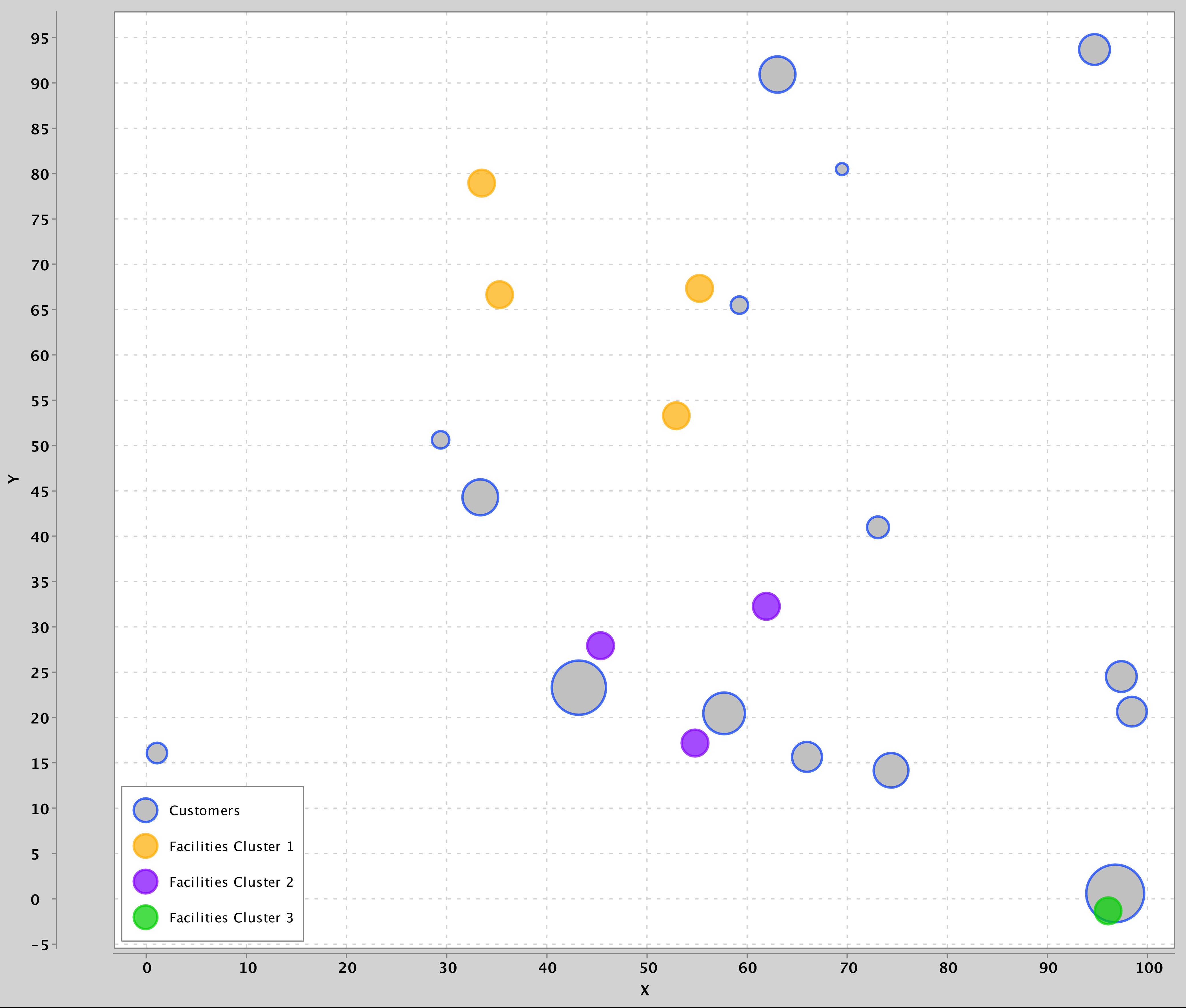}
  \caption{Example instance with eight facilities and fifteen customers. Facilities are clustered into three zones.
    The size of the customers represents the base mean demand $\mu_j$.}
  \label{fig:ex:instance}
\end{figure}

In this section we discuss the potential practical impact of accounting for decision-dependent uncertainty.
We do this by means of small illustrative examples. These facilitate illustration and comprehension of the underlying dynamics.
We use an instance with eight potential locations and fifteen customers. This is depicted in \Cref{fig:ex:instance}.
We assume locations are partitioned into three zones (clusters).

We consider demand type $\texttt{dt}=A$. The same analysis for demand types $B$, $C$, and $D$ is provided in \Cref{sec:Appendix2}.

We perform the following experiment. (i) We start by solving a simplified problem where the effect of location decisions on customers demand is neglected.
This is done by assuming that all facilities belong to the same zone (on the same instance), thus $\vert\mtc{Z}\vert =1$.
The solution to this problem is saved.
(ii) Next, we solve the original problem (i.e., where we account for decision-dependent uncertainty) but fix location decisions to those made in the simplified problem.
(iii) Finally, we solve the original problem without fixing location decisions.
This allows us to compare solutions from points (i) and (iii) and to assess the performance of the simplified solution from (i) into the actual problem where decisions have an impact on the uncertainty.

\begin{figure}
  \centering
  \begin{subfigure}[b]{0.45\textwidth}
    \centering
    \includegraphics[width=\textwidth]{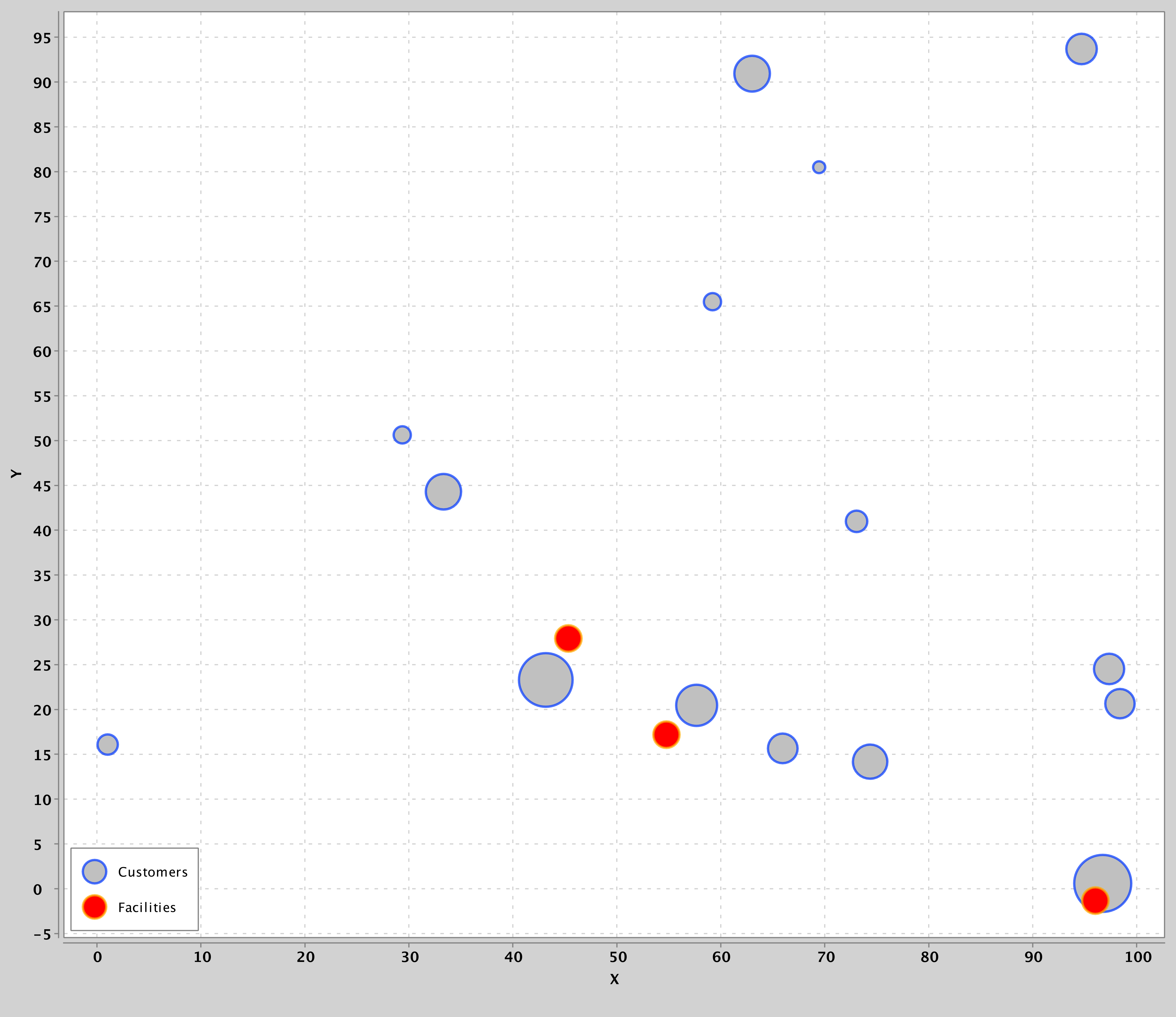}
    \caption{Solution that disregards decision-dependent uncertainty.}
    \label{fig:ex:solA1z}
  \end{subfigure}
  \hfill
  \begin{subfigure}[b]{0.45\textwidth}
    \centering
    \includegraphics[width=\textwidth]{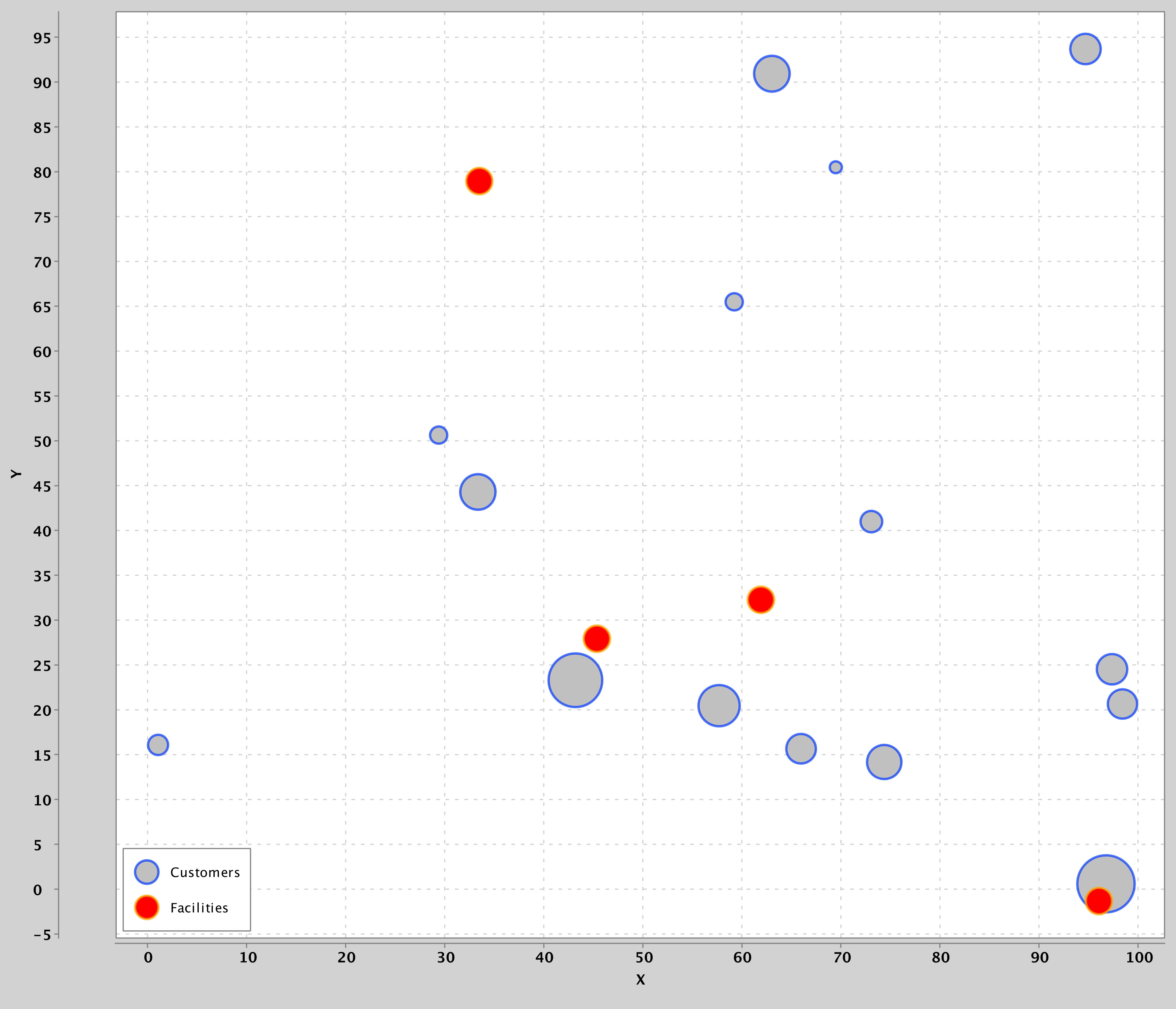}
    \caption{Solution which accounts for decision-dependent uncertainty.}
    \label{fig:ex:solA3z}
  \end{subfigure}
  \caption{Solutions to the instance depicted in \Cref{fig:ex:instance} with $\texttt{dt}=A$ obtained by disregarding (a) and considering (b) endogenous uncertainty.}
  \label{fig:ex:solutionsA}
\end{figure}

In \Cref{fig:ex:solutionsA} we report two solutions.
The solution in \Cref{fig:ex:solA1z} is obtained assuming that all facilities belong to the same geographical zone (see point (i) above), thus $\vert\mtc{Z}\vert=1$.
This gives us a facility location problem with fully exogenous uncertainty. We hereby effectively neglect the effect of decisions on the uncertainty.
The resulting demand distributions will have means $\hat{\mu}_{j}=\mu_j+\alpha^1\mu_j$ and variances $\hat{\sigma}_{j}=\sigma_j-\beta^1\sigma_j$ for all $j\in\mtc{J}$.
The solution in \Cref{fig:ex:solA3z} accounts for the endogenous effect of facility location decisions on the demand distribution by solving an instance of problem \eqref{eq:flp:compact} with $\vert\mtc{Z}\vert=3$ (see point (iii) above).

The main difference between the solutions in \Cref{fig:ex:solA1z} and \Cref{fig:ex:solA3z} is, essentially,
that the solution in \Cref{fig:ex:solA3z} opens an additional facility located in cluster 1.
This additional facility has a positive effect on the demand.
Particularly, the mean demand of the customers in the top-right region will increase (and the variance decrease) since their closest zone has at least one facility open.
In addition, the remaining customers will have open facilities in their $2$-nd or $3$-rd closest zone.
This will likewise increase their expected demand and decrease variance.

The difference between the solutions has an impact on expected profits.
The expected profit of the location decisions that account for endogenous uncertainty amounts to $272416.17$.
The expected profit of the solution to the simplified problem is, instead, $233538.71$.
This expected profit is calculated by fixing the solution in \Cref{fig:ex:solA1z} in the ``true'' problem with endogenous uncertainty, see point (ii).
In this simple example, accounting for the endogenous effect of facility locations on the uncertainties, generates a $16.4\%$ profit increase.

The impact of accounting for decision-dependent uncertainty is likely problem-specific. In the proposed problem is likely influenced by the $\alpha$ and $\beta$ parameters.
Therefore, we next assess how sensitive decisions are to the parameters $\alpha$ and $\beta$.
These determine how the mean and standard deviation of the demand distribution change as location decisions change.
To do so, we repeat the same experiment, but this time we use different coefficients $\alpha_{nj}$ and $\beta_{nj}$.
The base values considered in our experiments are $(\alpha_{nj},\beta_{nj})=(0.5,0.4)$. (Recall that in the formulas of $\mu_{jd}$ and $\sigma_{jd}$ they are raised to an exponent and decrease exponentially with the distance from the customer, see \Cref{sec:exp:inst}).
We now consider two additional scenarios where both $\alpha_{nj}$ and $\beta_{nj}$ decrease substantially with respect to the base values.
This mimics a situation where the impact of location decisions on the parameters of the distributions is considerably reduced.
In particular, we consider a scenario with $(\alpha_{nj},\beta_{nj})=(0.25,0.20)$ (hence the impact of location decisions is halved) and a scenario with $(\alpha_{nj},\beta_{nj})=(0.1,0.08)$ (hence the impact of location decisions is reduced to one fifth).
In these cases we observe that accounting for decision-dependent uncertainty generates a $6.6$\% and a $2.4$\% increase in profit, respectively, compared to the $16.4\%$ of the base case. This also suggests that, if precise estimates of the parameters $\alpha$ and $\beta$ cannot be obtained, one may consider developing Robust or Distributionally Robust extensions (see, e.g., \citet{BasAS21}) of the proposed method. These techniques are explicitly designed to plan against worst-case occurrences of parameters values and their distributions, respectively.

\section{Conclusions}\label{sec:conclusions}

This article presented a framework for modeling facility location problems characterized by decision-dependent uncertainty and an exact solution method.
The modeling framework allows specifying probability distributions that depend on the facilities open.
Particularly, given a partition of the facilities, the framework allows the distribution to change depending on e.g., which groups of facilities have been open.
It includes the limiting case where every facility location decision determines a different distribution.
The solution method extends the L-Shaped method by allowing distributions dependent on the first-stage solutions.
For this method, we provided distribution-specific optimality cuts and proved the cuts effective and safe.
Extensive tests on more than three thousand instances of a facility location problem with decision-dependent demand distribution provided empirical evidence of the efficency of the method.
We were able to solve to provable optimality instances with up to $2^{17}$ potential distributions. We have also shown that the performance of the method is affected by the coarsness of the partition of the facilities.

The strategy adopted in this article identifies a promising avenue of research for other combinatorial optimization problems characterized by decision-dependent uncertainty.
The central characteristic of the facility location problem treated in this article is that the set of feasible solution is countable and finate.
Distributions are enforced by subsets of these solutions, of which we naturally have a finite number. This generates a combinatorial number of potential distributions which makes enumeration inapplicable on all but trivial cases.
These characteristic can be identified in several other combinatorial optimization problems characterized by decisions such as selection, sequencing, matching.
Therefore, the strategy adopted in this article appears viable, with the necessary adjustments, in a much broader domain that that of location decisions.

\section*{Declarations \& Statements}

The procedure to obtain the datasets generated during and/or analysed during the current study are fully described in \Cref{sec:exp:inst}.

The author has no relevant financial or non-financial interests to disclose.

This research was partly supported by the \textit{Novo Nordisk Fonden} grant NNF24OC0089770.

\begin{appendices}
\crefalias{section}{appsec}
\section{Extensive Linearized Formulation}  \label{sec:Appendix1}
  In this section we report the extensive linearized formulation of problem \eqref{eq:flpsd:compact}. We introduce decision variables $w_{ij}^{sd}$ to represent the flow of product between facility $i$ and customer $j$ under realization $s\in\mtc{S}_d$ of distribution $d\in\mtc{D}$. 
  \begin{subequations}
    \label{eq:pmax:ext}
  \begin{align}
    \min &\sum_{i\in\mtc{I}}F_i x_i-\sum_{d\in\mtc{D}}\delta_d\sum_{s\in\mtc{S}_d}\pi_{sd}\left(\sum_{i\in\mtc{I}}\sum_{j\in\mtc{J}}R_{ij}w_{ij}^{ds}\right)\\
    s.t.~&\sum_{z\in\mtc{Z}_d}y_z-\sum_{z\in\mtc{Z}_d^C}y_z \geq |\mtc{Z}_d| \delta_d- |\mtc{Z}_d^C|(1-\delta_d)& \forall d \in\mtc{D}\\
         &\sum_{z\in\mtc{Z}_d}y_z-\sum_{z\in\mtc{Z}_d^C}y_z \leq \delta_d + |\mtc{Z}_d|-1& \forall d \in\mtc{D}\\
         &\sum_{i\in\mtc{I}_z}x_i\leq |\mtc{I}_z| y_z  &\forall z\in\mtc{Z}\\
         &\sum_{i\in\mtc{I}_z}x_i\geq y_z  &\forall z\in\mtc{Z}\\
         &\sum_{i\in\mtc{I}}w_{ij}^{sd}\leq \xi_j^{sd} & j\in\mtc{J},d\in\mtc{D},s\in\mtc{S}_d\\    
         &\sum_{j\in\mtc{J}}w_{ij}^{sd}\leq C_ix_i& i\in\mtc{I},d\in\mtc{D},s\in\mtc{S}_d\\    
         &w_{ij}^{sd}\geq 0&i\in\mtc{I},j\in\mtc{J} ,d\in\mtc{D},s\in\mtc{S}_d\\
         &x_i\in\{0,1\}&\forall i\in\mtc{I}\\
         &y_z\in\{0,1\}&\forall z\in\mtc{Z}\\
         &\delta_d\in\{0,1\}&\forall d \in\mtc{D}
  \end{align}
\end{subequations}
The bilinear products $\delta_dw_{ij}^{ds}$ in the objective function can be linearized by introducing new variables $h_{ij}^{sd}$ with the following additional constraints
\begin{equation*}
\begin{rcases}
  & h_{ij}^{sd}\leq w_{ij}^{sd}\\
  & h_{ij}^{sd}\leq \min\{\xi_j^{sd},C_i\}\delta_d \\    
  & w_{ij}^{sd}-h_{ij}^{sd}+\min\{\xi_j^{sd},C_i\}\delta_d\leq \min\{\xi_j^{sd},C_i\}\\
\end{rcases} i\in\mtc{I},j\in\mtc{J} ,d\in\mtc{D},s\in\mtc{S}_d
\end{equation*}

\section{Background Methodology}\label{app:ls}

A linear two-stage stochastic program can be formulated as follows.
\begin{equation}\label{eq:sp}
  \min_{x\in\mtc{X}}\{c^\top x + Q(x)\}
\end{equation}
where $\mtc{X}\subseteq \mathbb{R}^{n_1}$ is a polyhedral set, $Q(x)=\mathbb{E}_{\pr{P}}\left[Q(x,\rv{\xi})\right]$ and, for a given $x$ and a given realization $\xi$ of $\rv{\xi}$
$$Q(x,\xi)=\min_{y\in \mtc{Y}}\{q^\top y\vert Wy=h-Tx\}$$
represents the cost of second-stage decisions. Here $\mtc{Y}\subseteq \mathbb{R}^{n_2}$ is a polyhedral set and $\rv{\xi}\sim \pr{P}$, independent of $x$, collects the random elements of the second-stage problem $Q(x,\rv{\xi})$, that is $\rv{q},\rv{W},\rv{T},\rv{h}$.

The L-Shaped method solves problems of type \eqref{eq:sp} under specific assumptions, see \cite{VanW69}.
In particular, it requires that $\rv{\xi}$ has a discrete distribution independent of $x$.
Let $\mtc{S}$ be the collection of possible scenarios of $\rv{\xi}$.
That is, we have realization $\xi_s$ and probabilities $\pi_s$, $s\in\mtc{S}$. Then
$Q(x)=\sum_{s\in\mtc{S}}\pi_sQ(x,\xi_s)$.

\Cref{alg:app:ls} sketches the main steps of the method.
At each iteration, the method works with the following relaxation of \eqref{eq:sp}
\begin{equation}\label{eq:ls:mp}
  \min_{x\in\mtc{X},\mu\in \mathbb{R}}\{c^\top x + \mu\}
\end{equation}
Problem \eqref{eq:ls:mp} is referred to as the \textit{Relaxed Master Problem} (RMP). The variable $\mu$ approximates $Q(x)$ from below.

\begin{algorithm}
\caption{L-Shaped Method}\label{alg:app:ls}
\begin{algorithmic}[1]
  \State \texttt{solved}$\gets$\texttt{false}, $i\gets 0$.
  \While{not \texttt{solved}}
    \State Solve RMP. Let $(x^i,\mu^i)$ be its optimal solution.
    \State Solve $Q(x^i,\xi_s)$ for all $s\in\mtc{S}$.
    \If{$Q(x^i,\xi_s)$ is infeasible for some $s\in\mtc{S}$}
      \State Add a \textit{feasibility cut} to RMP.
    \Else
      \State Compute $\mathbb{E}_{\mathsf{P}}[Q(x^i,\rv{\xi})]$.
      \If{$\mu^i<\mathbb{E}_{\mathsf{P}}[Q(x^i,\rv{\xi})]$ }
        \State Add an \textit{optimality cut} to RMP.
      \Else   
        \State \texttt{solved}$\gets$\texttt{true}.
      \EndIf
    \EndIf
    \State $i\gets i+1$.
  \EndWhile
\end{algorithmic}
\end{algorithm}

At a given iteration $i$, one solves \eqref{eq:ls:mp} and obtains a solution $(x^i,\mu^i)$.
Three cases may materialize:
\begin{enumerate}
\item The second-stage problem $Q(x^i,\xi_s)$ is infeasible for some $s\in\mtc{S}$.
  In this case, the proposed first-stage solution $x^i$ renders some second-stage problem infeasible and is therefore not acceptable.
  The method describes how to generate a linear inequality in $x$, called a \textit{feasibility cut}.
  This is generated using an extreme ray to the dual of linear program $Q(x^i,\xi_s)$.
  This inequality is violated by $x^i$ and is added to \eqref{eq:ls:mp} starting from the next iteration.
\item The second-stage problem $Q(x^i,\xi_s)$ is feasible for all $s\in\mtc{S}$ but $\mu^i<Q(x^i)$.
  In this case, $\mu^i$ does not accurately estimate the expected cost of $x^i$.
  It must, therefore, be corrected. The method describes how to generate a linear inequality in $x$, called an \textit{optimality cut}, that is violated by $(x^i,\mu^i)$.
  This inequality is based on strong duality arguments regarding the linear programs $Q(x^i,\xi_s)$, $s\in\mtc{S}$.
  The optimality cut is added to \eqref{eq:ls:mp} from the next iteration to prevent producing solution $(x^i,\mu^i)$ again.
\item The second-stage problem $Q(x^i,\xi_s)$ is feasible for all $s\in\mtc{S}$ and $\mu^i=Q(x^i)$. Then $(x^i,\mu^i)$ is optimal for \eqref{eq:ls:mp}.
\end{enumerate}
The method is guaranteed to converge to an optimal solution within a finite number iterations. 
Several extensions and refinements of the method exist.

The extension developed by \cite{PanH25} allows solving problems where the distribution of $\rv{\xi}$ depends on $x$.
In particular, the authors assume that the relationship between $x$ and the distribution $\pr{P}$ of $\rv{\xi}$ is piecewise-constant on $\mtc{X}$.
That is, there exists a finite partition $(\mtc{X}_d)_{d\in\mtc{D}}$ of $\mtc{X}$, where $\mtc{D}$ is a finite discrete set, such that $\pr{P}$ is constant (i.e., the same distribution applies) for all choices of $x$ in a given subset  $\mtc{X}_d$ of the partition. In other words, decisions $x\in\mtc{X}_d$ enforce distribution $\pr{P}_d$, for $d\in\mtc{D}$.
Assume $\pr{P}_d$ is discrete for all $d\in\mtc{D}$, with scenarios $s\in \mtc{S}_d$, and realizations $\xi_{sd}$.

The revised method is sketched in \Cref{alg:app:lsddu}. The L-Shaped method is adjusted as follows. 
Upon obtaining a solution $(x^i,\mu^i)$ to RMP, the method identifies the distribution $\pr{P}_{d^i}$, $d^i\in\mtc{D}$ enforced by $x^i$.
From that point on, and for the rest of iteration $i$, the method works only with that distribution. In particular,
\begin{itemize}
\item It verifies whether the second-stage problems are feasible \textit{only} under the possible realizations of $\rv{\xi}$ in distribution $d^i$.
  If $Q(x^i,\xi_{sd^i})$ is infeasible for some $s\in\mtc{S}_{d^i}$, the method adds a \textit{distribution-specific} feasibility cut. This cut is of type
  $$x\in \mtc{X}_{d^i}\implies f(x)\leq 0$$
  That is, it applies only when the decision $x$ enforces distribution $d^i$ and is redundant otherwise. Here, $f(x)$ is an affine function of $x$. 
\item It verifies whether $\mu^i< \mathbb{E}_{\mathsf{P}_{d^i}}[Q(x^i,\rv{\xi})] $. That is, $\mu$ should now approximate the conditional (on $x$) expectation of $Q(x,\rv{\xi})$.
  If the condition is verified it adds a \textit{distribution-specific} optimality cut. This cut is of type
  $$x\in \mtc{X}_{d^i}\implies o(x)\leq \mu$$
  That is, the cut applies only when the decision $x$ enforces distribution $d^i$ and is redundant otherwise. Here, $o(x)$ is an affine function of $x$.
\end{itemize}
The method terminates when $\mu$ correctly estimates the conditional expectation of the second-stage cost.

\begin{algorithm}
\caption{L-Shaped Method for stochastic programs with decision-dependent uncertainty}\label{alg:app:lsddu}
\begin{algorithmic}[1]
  \State \texttt{solved}$\gets$\texttt{false}, $i\gets 0$.
  \While{not \texttt{solved}}
    \State Solve RMP. Let $(x^i,\mu^i)$ be its optimal solution.
    \State Identify the distribution $d^i\in\mtc{D}$ such that $x^i \in \mtc{X}_{d^i}$.
    \State Solve $Q(x^i,\xi_{sd^i})$ for all $s\in\mtc{S}_{d^i}$.
    \If{$Q(x^i,\xi_{sd^i})$ is infeasible for some $s\in\mtc{S}_{d^i}$}
      \State Add a \textit{distribution-specific} feasibility cut to RMP.
    \Else
      \State Compute $\mathbb{E}_{\mathsf{P}_{d^i}}[Q(x^i,\rv{\xi})]$.
      \If{$\mu^i<\mathbb{E}_{\mathsf{P}_{d^i}}[Q(x^i,\rv{\xi})]$ }
        \State Add an \textit{distribution-specific} optimality cut to RMP.
      \Else   
        \State \texttt{solved}$\gets$\texttt{true}.
      \EndIf
    \EndIf
    \State $i\gets i+1$.
  \EndWhile
\end{algorithmic}
\end{algorithm}

The condition that there exists a finite partition $\mtc{X}=\cup_{d\in \mtc{D}}\mtc{X}_d$, with $\mtc{X}_{d_1}\cap\mtc{X}_{d_2}=\emptyset$ for all $d_1,d_2\in \mtc{D}$ naturally holds in combinatorial optimization problems. In particular, as a special case, it holds when the decision variables may take on values from a finite set, and there are only a finite (though very large) number of possible combinations of these values. This is the case, for example, of the facility location problem studied in this paper.

\section{Additional Results on Potential Impact}\label{sec:Appendix2}

We provide additional evidence on the impact of decision-dependent uncertainty on profits under demand types $B$, $C$ and $D$.

\begin{figure}
  \centering
  \begin{subfigure}[b]{0.45\textwidth}
    \centering
    \includegraphics[width=\textwidth]{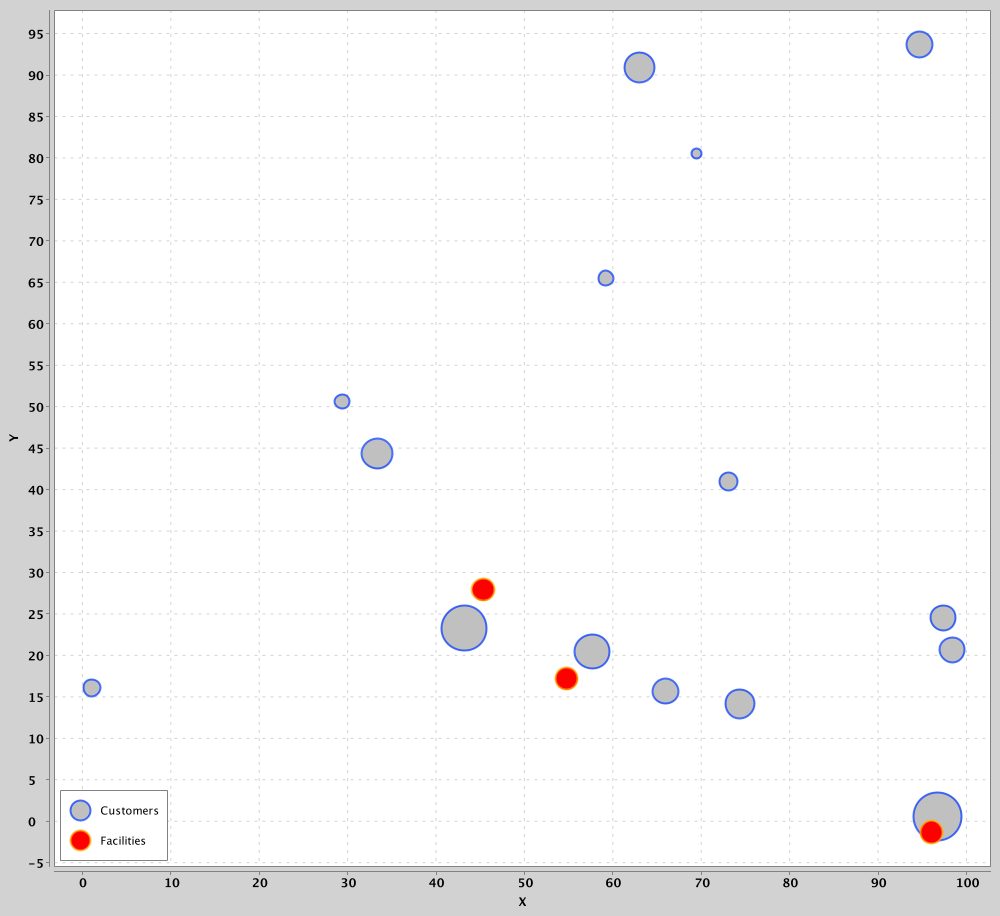}
    \caption{Solution that disregards decision-dependent uncertainty.}
    \label{fig:ex:solB1z}
  \end{subfigure}
  \hfill
  \begin{subfigure}[b]{0.45\textwidth}
    \centering
    \includegraphics[width=\textwidth]{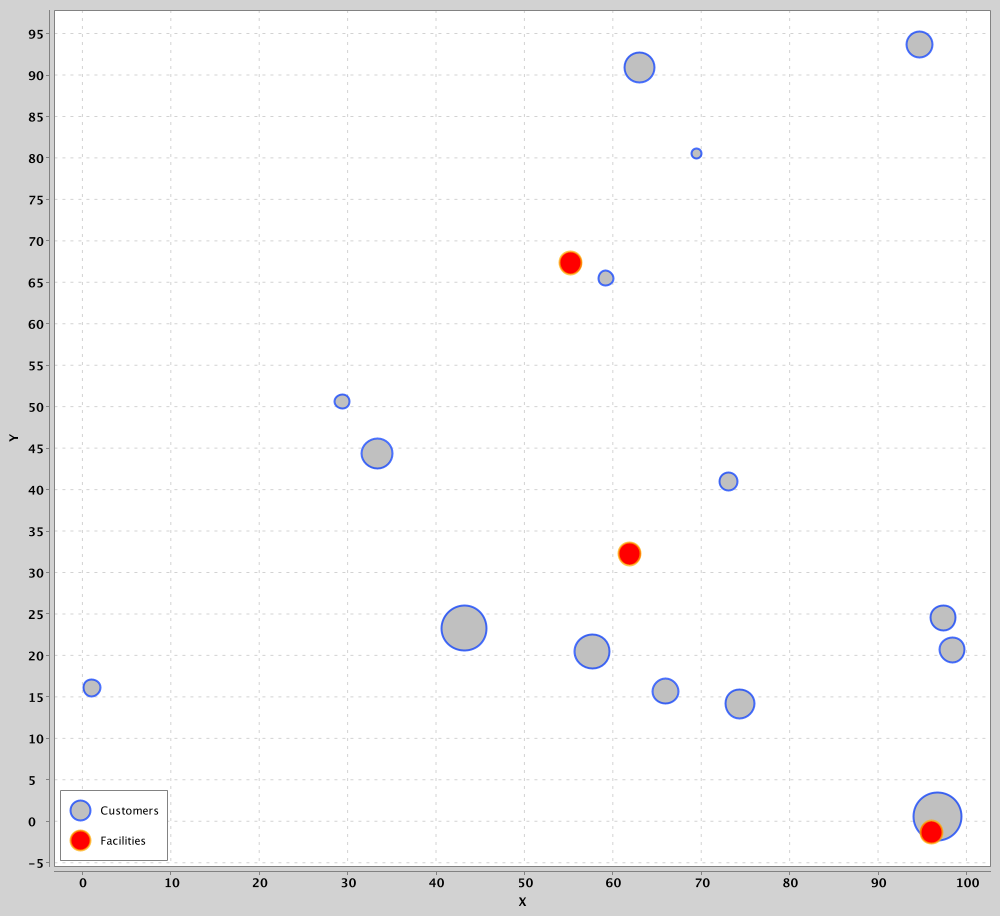}
    \caption{Solution which accounts for decision-dependent uncertainty.}
    \label{fig:ex:solB3z}
  \end{subfigure}
  \caption{Solutions to the instance depicted in \Cref{fig:ex:instance} with $\texttt{dt}=B$ obtained by disregarding (a) and considering (b) endogenous uncertainty.}
  \label{fig:ex:solutionsB}
\end{figure}

We start by using demand type $\texttt{dt}=B$. We remind the reader that, in this case, the demand increases (and variance decreases) if and only if there are open facilities in the nearest zone. Further zones have no effect. The solutions are reported in \Cref{fig:ex:solutionsB}.
We notice that both solutions suggest opening exactly three facilities. The main difference is in the zones they activate.
The solution that disregards decision-dependent uncertainty opens two facilities in cluster 2 and one in cluster 3, see \Cref{fig:ex:solB1z}. This solution does not account for the fact that the demand of the customers in the top-right corner would increase if their closest zone (cluster 1) was active.
This effect is instead accounted for when including decision-dependent uncertainty, see \Cref{fig:ex:solB3z}.
This slight difference in the solutions has a significant effect on expected profits.
The solution that accounts for decision-dependent uncertainty generates an $11.2\%$ higher profit.

\begin{figure}
  \centering
  \begin{subfigure}[b]{0.45\textwidth}
    \centering
    \includegraphics[width=\textwidth]{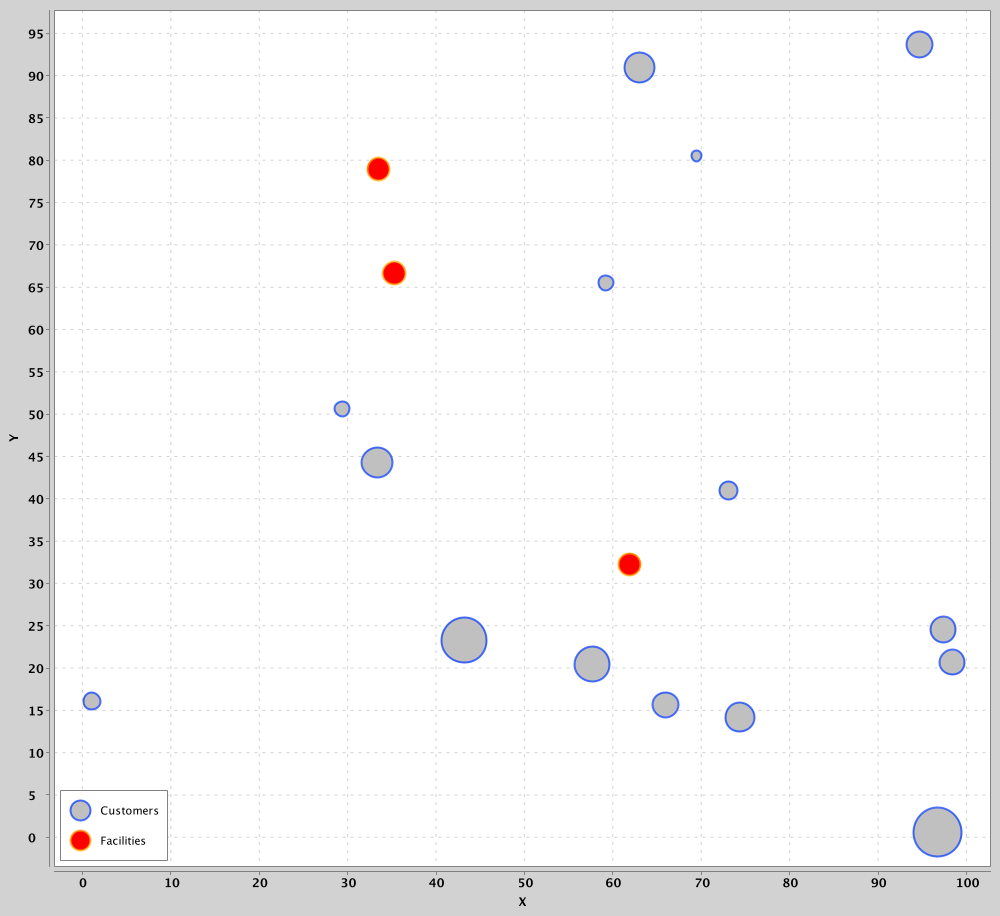}
    \caption{Solution that disregards decision-dependent uncertainty.}
    \label{fig:ex:solC1z}
  \end{subfigure}
  \hfill
  \begin{subfigure}[b]{0.45\textwidth}
    \centering
    \includegraphics[width=\textwidth]{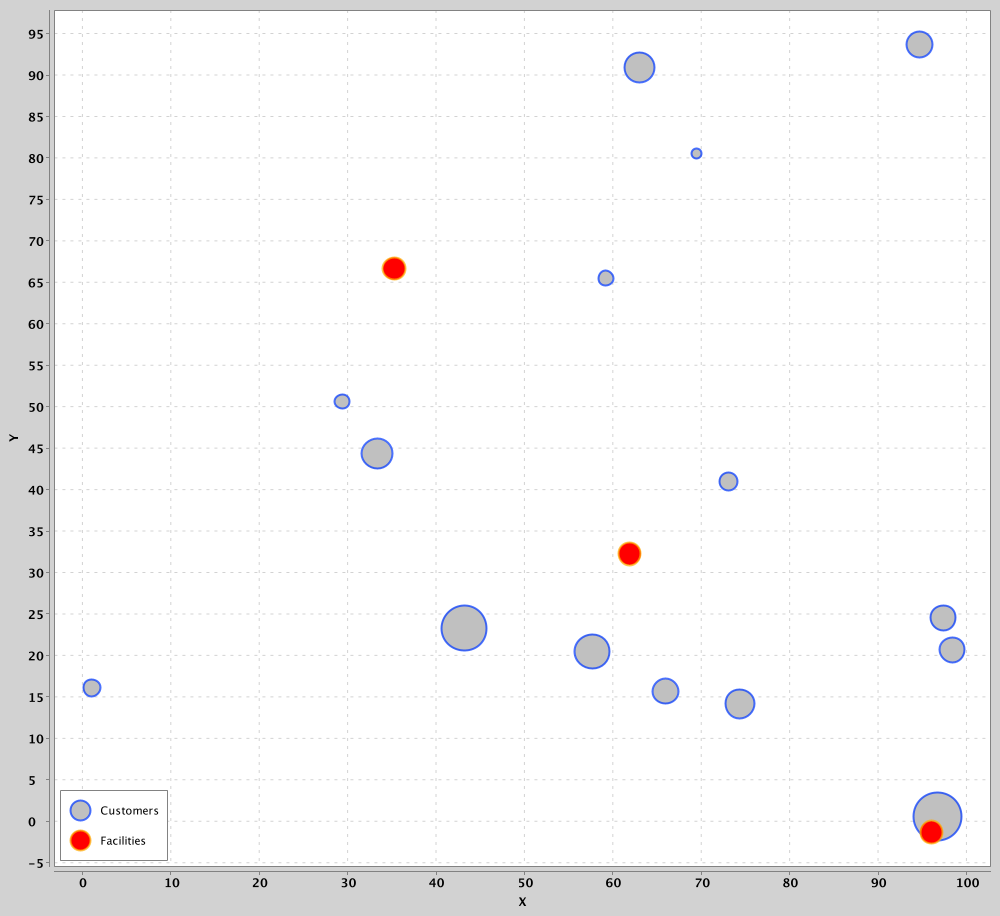}
    \caption{Solution which accounts for decision-dependent uncertainty.}
    \label{fig:ex:solC3z}
  \end{subfigure}
  \caption{Solutions to the instance depicted in \Cref{fig:ex:instance} with $\texttt{dt}=C$ obtained by disregarding (a) and considering (b) endogenous uncertainty.}
  \label{fig:ex:solutionsC}
\end{figure}

We then perform the same experiment using demand type $\texttt{dt}=C$. The solutions are reported in \Cref{fig:ex:solutionsC}.
We remind the reader that under $\texttt{dt}=C$ the mean of the demand increases (and the variance decreases) if and only if
the closest open zone is the $n$-furthest, by a percentage that decreases as $n$ increases. In \Cref{fig:ex:solutionsC} we observe that both the solution that disregards decision-dependent uncertainty (\Cref{fig:ex:solC1z}) and the solution that account for that (\Cref{fig:ex:solC3z}) indicate opening three out of the eight facilities. The only difference in the solution is in the choice of locations.
The solution in \Cref{fig:ex:solC1z} opens two facilities in cluster 1 and one in cluster 2. The solution in \Cref{fig:ex:solC3z} opens one facility in each cluster.
Cluster 3 is the closest to the customer in the bottom-right corner. It appears that the problem is able to identify that opening the closest facility to that customer would increase its demand, and deems it profitable.
Indeed, the expected profit of the solution that accounts for decision-dependent uncertainty is $4.0\%$ higher.

\begin{figure}
  \centering
  \begin{subfigure}[b]{0.45\textwidth}
    \centering
    \includegraphics[width=\textwidth]{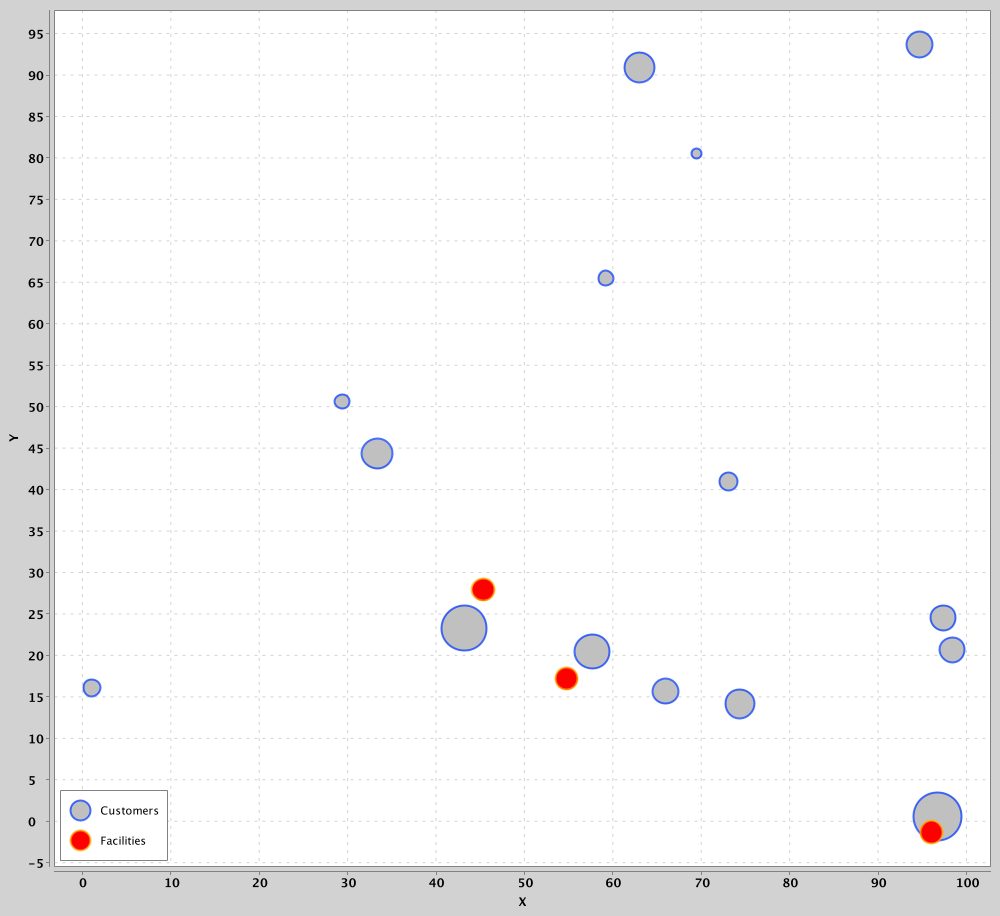}
    \caption{Solution that disregards decision-dependent uncertainty.}
    \label{fig:ex:solD1z}
  \end{subfigure}
  \hfill
  \begin{subfigure}[b]{0.45\textwidth}
    \centering
    \includegraphics[width=\textwidth]{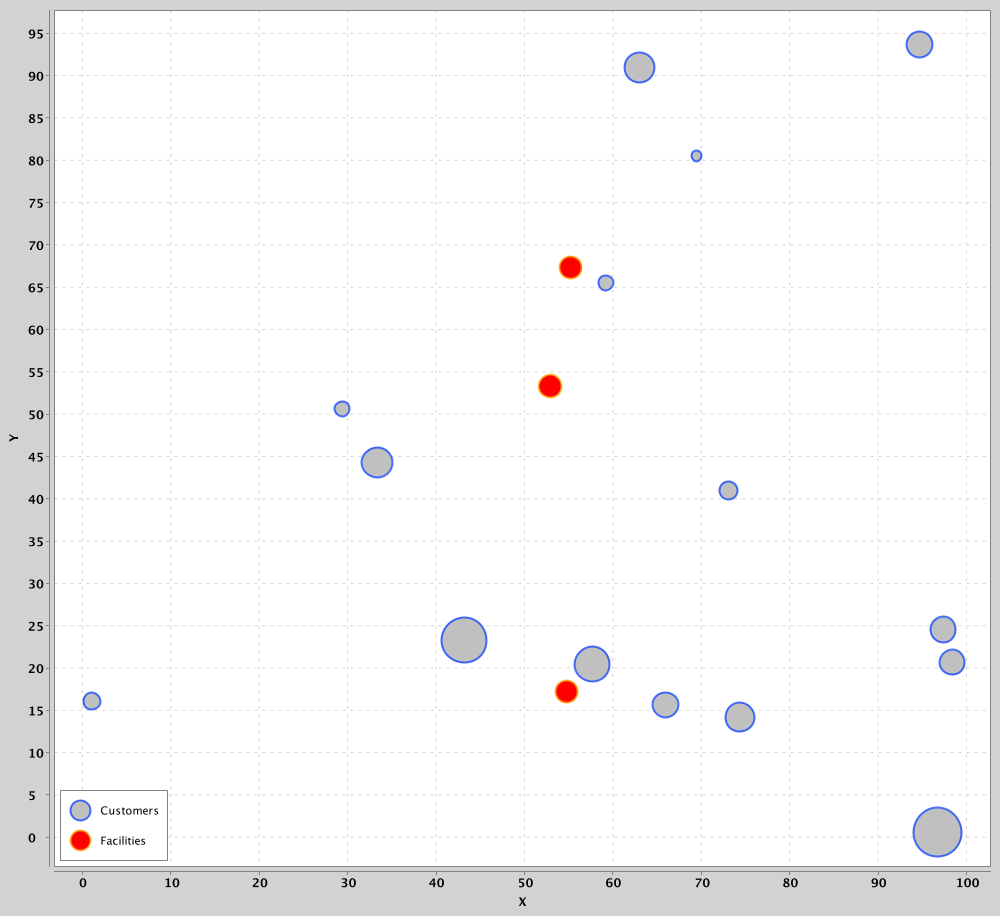}
    \caption{Solution which accounts for decision-dependent uncertainty.}
    \label{fig:ex:solD3z}
  \end{subfigure}
  \caption{Solutions to the instance depicted in \Cref{fig:ex:instance} with $\texttt{dt}=D$ obtained by disregarding (a) and considering (b) endogenous uncertainty.}
  \label{fig:ex:solutionsD}
\end{figure}

Finally, we consider demand type $\texttt{dt}=D$. The solutions are reported in \Cref{fig:ex:solutionsD}.
We remind the reader that, under $\texttt{dt}=D$, the mean of the demand of a customer increases (and the variance decreases) if the closest zone is open and decreases (the variance increases) if zones further away are open.
That is, the opening of facilities far away from the customer contributes negatively to the demand of the customer.
In \Cref{fig:ex:solutionsD} we observe that both the solution that disregards decision-dependent uncertainty (\Cref{fig:ex:solD1z}) and the solution that accounts for that (\Cref{fig:ex:solD3z}) indicate opening three out of the eight facilities.
In addition, both solutions indicate covering two of the three zones. The only difference is in the choice of open zones.
The solution in \Cref{fig:ex:solD1z} opens facilities in clusters 2 and 3. The solution in \Cref{fig:ex:solD3z} opens one facility in clusters 1 and 2.
A plausible explanation is the following. There are more customers close to cluster 1 than to cluster 3. Opening facilities in cluster 1 will decrease the demand of the customers close to cluster 3. Opening facilities in cluster 3 will decrease the demand of the customers close to cluster 1.
Since the customers close to cluster 3 are fewer, the solution suggests opening facilities in cluster 1. 
The expected profit of the solution that accounts for decision-dependent uncertainty is $3.8\%$ higher.

\end{appendices}

\bibliographystyle{plainnat} 
\bibliography{gpbd}

\end{document}